\documentclass{amsart}
\usepackage{amsmath,amssymb,amsthm,amscd,mathtools}
\usepackage[dvipsnames]{xcolor}
\usepackage{tikz-cd}
\usepackage{enumitem}
\usepackage{array,booktabs,calc}
\usepackage[textwidth=6.1in,textheight=8.9in,centering]{geometry}
\usepackage{microtype}
\usepackage{hyperref}
\hypersetup{
colorlinks=true,
pdftitle={\detokenize{The RO(C_{p^n})-graded Tambara functor of perfectoid THH}},
pdfauthor={Dongming Zhang},
linkcolor=blue,
anchorcolor=blue,
citecolor=blue
}
\usepackage[nameinlink]{cleveref}

\newtheorem{theorem}{Theorem}[section]

\newtheorem{proposition}[theorem]{Proposition}
\newtheorem{lemma}[theorem]{Lemma}

\theoremstyle{definition}
\newtheorem{definition}[theorem]{Definition}

\theoremstyle{remark}

\makeatletter
\let\c@equation\c@theorem
\makeatother
\numberwithin{equation}{section}

\newcommand{\HZ}{H\underline{\mathbb Z}}
\newcommand{\res}{\operatorname{res}}
\newcommand{\Res}{\operatorname{Res}}
\newcommand{\tr}{\operatorname{tr}}
\newcommand{\Ind}{\operatorname{Ind}}

\newcommand{\coker}{\operatorname{coker}}
\newcommand{\Sp}{\operatorname{Sp}}

\newcommand{\holim}{\operatorname*{holim}}

\newcommand{\THH}{\operatorname{THH}}
\newcommand{\TF}{\operatorname{TF}}
\newcommand{\TR}{\operatorname{TR}}
\newcommand{\RO}{RO}
\newcommand{\Z}{\mathbb Z}
\newcommand{\F}{\mathbb F}
\newcommand{\Map}{\operatorname{Map}}
\begin{document}
\author{Dongming Zhang}
\title[The $RO(C_{p^n})$-graded Tambara functor of perfectoid THH]{The $RO(C_{p^n})$-graded Tambara functor of perfectoid THH}
\begin{abstract}
Let $R$ be a perfect $\mathbb F_p$-algebra or a nonzero
$p$-torsionfree perfectoid ring and $p$ be an arbitrary prime number. We compute the $RO(C_{p^n})$-graded
Tambara functor of $\operatorname{THH}(R;\mathbb Z_p)$ for every $n\ge0$, including all homotopy groups, products, restrictions,
transfers, Weyl actions and norms.
We also compute genuine cyclotomic restriction. 
\end{abstract}
\maketitle
\setcounter{tocdepth}{1}
\tableofcontents

\section{Introduction}\label{sec:introduction}
Let $p$ be a prime and let $R$ be a perfect $\F_p$-algebra or a nonzero
$p$-torsionfree perfectoid ring. Put
\[
X_R=\THH(R;\Z_p),\qquad
\pi_V^H X_R=[S^V,X_R]^H\quad(V\in RO(H)).
\]
We use the genuine commutative cyclotomic model of $\THH(HR)$
from \cite[Theorems~4.6--4.7]{NormTHH} and
\cite[Theorems~C--E]{BMCyclotomic}, with
\[
X_R=\holim_m\THH(HR)/p^m.
\]
The quotients are taken in spectra.
On geometric fixed points, a norm is a finite smash product over double cosets
\cite[Theorem~2.35]{NormTHH}. Thus norms preserve genuine mod-$p$ equivalences, and
\cite[Theorem~6.1]{BHNorm} implies the existence of norms on $X_R$.
On finite cyclic $p$-groups, isotropy separation and
\cite[Lemma~2.11(b)]{AN} show that geometric $C_p$-fixed points
commute with this uniformly connective tower.
Thus
\[
\varrho_p^*\Phi^{C_p}X_R\simeq X_R,
\qquad
\pi_k\Phi^{C_{p^i}}X_R=0
\quad(i\ge0,\ k<0).
\]
Thus $X_R$ admits a genuine cyclotomic restriction.

We can now state our main result.

\begin{theorem}\label{thm:main}
Let $p$ be a prime, and let $R$ be either a perfect
$\F_p$-algebra or a nonzero $p$-torsionfree perfectoid
ring. For every $G=C_{p^n}$, $n\ge0$, the $RO(G)$-graded
Tambara functor and genuine cyclotomic restriction of $X_R$ are given as follows.
\begin{enumerate}[label=\textup{(\roman*)}]
\item For nonzero $p$-torsionfree $R$, the homotopy groups $\pi_V^H X_R$ for every
$H\le G$ and $V\in RO(H)$ are given by
Proposition~\ref{thm:uniform}.
Products, restrictions, and transfers are given by
Proposition~\ref{prop:products-maps}, and Weyl actions
by Proposition~\ref{prop:weyl-action}.

Norms are given by
Propositions~\ref{prop:vertex-norms},
\ref{lem:period-repair-norm}, and
\ref{prop:edge-norms}.
The extension to arbitrary representation degrees
is given by Proposition~\ref{lem:ordered-sign-units}.
Genuine cyclotomic restriction is given by
Proposition~\ref{prop:cyclotomic-operations}.

\item For a perfect $\F_p$-algebra $R$, Theorem~\ref{thm:cp-result}
gives the graded Tambara functor by Euler localization of $\HZ$ and
extension of coefficients to $W(R)$, and gives genuine cyclotomic
restriction in every representation degree.
\end{enumerate}
\end{theorem}
The main strategy we use is Tate compression introduced in \cite[Theorem~4.7]{ZhangHZ}, so we write the homotopy groups as a kernel and a
cokernel.

Our calculation is related to several areas, including topological cyclic homology,
algebraic $K$-theory, and $p$-adic Hodge theory.

Our calculation applies to topological cyclic
homology. Write
\[
\TR^{n+1}(R;\Z_p)=X_R^{C_{p^n}}.
\]
For $n\ge1$, there are Frobenius, Verschiebung, and cyclotomic
restriction maps
\[
F,R:\TR^{n+1}(R;\Z_p)\longrightarrow
\TR^n(R;\Z_p),
\qquad
V:\TR^n(R;\Z_p)\longrightarrow
\TR^{n+1}(R;\Z_p).
\]
On homotopy groups, these maps are given by
$\res$, $R_{\mathrm{cyc}}$, and $\tr$, respectively.
For perfectoid rings, see \cite[Proposition~7.21]{MathewTR} for the
finite-level integer-graded $\TR$ rings and \cite[Theorem~3.7]{Andriopoulos}
for the formulas for $R$ and $F$.
For odd $p$ and $p$-torsionfree perfectoid $R$ containing a compatible
system of primitive $p$-power roots of unity, Davis and
Patchkoria~\cite[Section~8]{DavisPatchkoria} identify
$\pi_2\TR^n(R;\Z_p)$ with the Tate module $T_p(W_n\Omega_R^1)$.

For a perfect $\F_p$-algebra $R$, we have
\[
\pi_0\TR^{n+1}(R;\Z_p)\cong W_{n+1}(R),
\]
where $W_{n+1}(R)$ is the ring of Witt vectors of
length $n+1$~\cite[Theorem~3.3]{HM97}.
Our formulas extend this calculation to all real representation
degrees. In degree zero our norm is the Witt vector
norm of Angeltveit~\cite[Theorem~4.1]{WittNorm}.
For related degree-zero calculations using Witt vectors for Green
functors and Witt vectors with coefficients, see
\cite{BGHLWitt,DKNPWitt}.
For shifts by actual complex representations over perfect fields,
see \cite[Proposition~9.1]{HM97}.
For $R=\F_p$, the $\TR$-groups in complex and real circle representation
degrees were computed in \cite[Theorem~1.1]{GerhardtTR} and
\cite[Section~7]{AGTR}, respectively.

The spectrum $\THH(R;\Z_p)$ also appears in algebraic
$K$-theory. The cyclotomic trace~\cite{BHM} gives a map
\[
K(R)^\wedge_p\longrightarrow\operatorname{TC}(R;\Z_p).
\]
Topological cyclic homology is constructed from the
cyclotomic structure on $\THH(R;\Z_p)$
\cite{NS,BMCyclotomicOriginal}.
Our result gives the homotopy groups and the maps
in its finite cyclic fixed-point tower.
Hesselholt and Madsen~\cite{HMPolytopes} used representation
suspensions of $\THH(k)$ to compute $K_*(k[x]/(x^e),(x))$ for
a perfect field $k$ of positive characteristic.

For $p$-torsionfree perfectoid rings, our calculation
is also related to $p$-adic Hodge theory and
prismatic cohomology.
Let
\[
A=W(R^\flat),\qquad
I_0=\ker(\theta:A\to R),\qquad
I_j=\varphi^j(I_0).
\]
These ideals come from the perfect prism associated
to $R$~\cite[Theorem~3.10]{BSPrisms}. The Breuil--Kisin twist $B$ satisfies
\[
\pi_2\operatorname{TP}(R;\Z_p)\cong B,
\qquad
\pi_2\operatorname{TC}^-(R;\Z_p)\cong I_0B.
\]
For these isomorphisms, the canonical map, and linearized cyclotomic
Frobenius, see \cite[Proposition~6.5 and Remark~6.6]{BMS}
and Subsection~\ref{subsec:period-line}.
Proposition~\ref{thm:uniform} expresses the groups
$\pi_V^{C_{p^n}}X_R$ using $I_j$ and $B$.

For other historical works,
Hesselholt~\cite[Theorem~B]{HesCurves} described the integer-graded
homotopy groups of finite-level $\TR$ in terms of de Rham--Witt forms
and a polynomial generator in degree two for smooth algebras over a perfect field of characteristic $p$.
For perfectoid rings $R$, Sulyma studied the regular slice filtration
of $\THH(R;\Z_p)$~\cite{SulymaSlice} and the $RO(\mathbb T)$-graded
$\TF$ of $R$~\cite{SulymaTF}. He formulated the homotopy orbits to
$\TR$ spectral sequence as a spectral sequence of Mackey functors
and discussed the $RO(\mathbb T)$-graded Tambara structure of
$\THH(R;\Z_p)$~\cite[Sections~1 and~3.2]{SulymaTF}.

Let $A$ be a $\delta$-ring with a spherical lift $\mathbb S_A$.
For relatively quasiregular semiperfectoid $A$-algebras $S$,
Sulyma~\cite[Theorem~1.9]{SulymaNew} computed the representation-graded
groups of $\operatorname{TC}^-(S\mid\mathbb S_A)$ and
$\operatorname{TP}(S\mid\mathbb S_A)$.
For a polyrepresentation degree $\alpha$ in the contiguous range, he also computed
the form-degree $j$ contributions to finite-level relative $\TR$ in total degrees
$\alpha+j$ and $\alpha-1+j$ for $p$-completely
smooth algebras over $p$-torsionfree perfectoid rings and for smooth
algebras over perfect $\F_p$-algebras~\cite[Theorem~4.12]{SulymaNew}.
His norm formulas for de Rham--Witt forms apply in all form degrees
at odd primes and in even form degrees at $p=2$%
~\cite[Remark~1.20 and Theorem~5.7]{SulymaNew}.
He deferred a general treatment of Breuil--Kisin and prismatic norms
for the homotopy fixed points and Tate spectra of relative
$\THH$~\cite[Section~1.3.3]{SulymaNew}.
For the rings $R$ in Theorem~\ref{thm:main}, we determine the genuine
norms of $\THH(R;\Z_p)$ for every inclusion $H\le G=C_{p^n}$ and
every degree $V\in RO(H)$, including degrees with an odd coefficient
of the sign representation at $p=2$.

In characteristic $p$, Krause and Nikolaus~\cite[Proposition~2.1]{KrauseNikolaus}
identified the underlying spectrum of $\THH(k)$ with
$k\otimes_{\F_p}\THH(\F_p)$ for every perfect $\F_p$-algebra $k$.
At odd primes, Ayala, Mazel-Gee, and Rozenblyum~\cite[Corollary~10.8]{AMRDerived}
identified $\THH(\F_p)$ with $\Phi^{C_p}\HZ$ as genuine spectra on
the tower of cyclic $p$-groups. For the multiplicative comparison, see
Mao~\cite[Remark~7.6]{MaoDecomplete}.

Hu, Kriz, and Somberg~\cite[Theorem~2]{HKSMotivic} described the
equivariant motivic filtration on $\THH$ of polynomial $\F_p$-algebras.
Burghardt--Hu--Kriz--Kriz--Somberg \cite{Kriz}
gave chain-level models for equivariant $\THH$
in characteristic $p$.
Our characteristic-$p$ comparison extends their
module comparison to preserve norms.

Bhatt--Morrow--Scholze \cite[Theorem~1.6]{BMS} computed the
integer-graded homotopy groups of $\THH(-;\Z_p)$,
$\operatorname{TC}^-(-;\Z_p)$, and $\operatorname{TP}(-;\Z_p)$
for perfectoid rings.
For odd $p$ and $R=\mathcal O_{\mathbb C_p}$, the identification of
$\pi_0\operatorname{TC}^-(R;\Z_p)$ with $W(R^\flat)$ goes back
to Hesselholt~\cite{HesLocal}.
Mao~\cite{MaoThom} gave a Thom spectrum proof of B\"okstedt
periodicity for perfectoid rings.
Riggenbach \cite{Riggenbach} compared circle
homotopy fixed points with finite-subgroup
homotopy fixed points.
Bhatt--Lurie \cite{BhattLurie} developed the
divided conormal construction, which is used for the
Breuil--Kisin twist.

For real $\THH$ of perfectoid rings and for equivariant Witt complexes
in twisted $\THH$, see \cite{HornbostelPark} and \cite{BGKPY}, respectively.
Other work on cyclotomic spectra,
norms, and Tambara functors includes
\cite{AN,NS,NormTHH,BMCyclotomic,BMYRelative,
BHNorm,Ullman,Balderrama,HKSZOperations}.

Our conventions are as follows.
We work in complete universes. Fixed points, Tate constructions,
limits and smash products are derived.
For $V\in RO(G)$, $S^V$ is its representation sphere.
An ordered model of $V$ fixes the order of its irreducible summands. For $C_N$, fix
$\gamma_N=\exp(2\pi i/N)$; for $C_M\le C_N$, we have
$\gamma_M=\gamma_N^{N/M}$.

\begin{definition}\label{def:base-class}
Write $\mathbb T=S^1$ and let $C_{p^n}$ be its subgroup of
$p^n$th roots of unity. Let
\[
\varrho_p:\mathbb T\xrightarrow{\sim}\mathbb T/C_p
\]
be the inverse of $[z]\mapsto z^p$.
For $H\le G$, a degree $V\in RO(G)$ is \emph{natural} at $G/H$
if $\dim V^H\in2\mathbb Z$ and \emph{shifted} otherwise.
We omit the orbit when $H=G$.
\end{definition}

The $\Res,\Ind$ act on representations;
$\res,\tr,c_g,N$ denote restriction, transfer, conjugation and norm.
For degrees in $RO(G)$, $c_g$ also denotes the Weyl action.
We write $R_{\mathrm{cyc}}$ for cyclotomic restriction and
$F(-,-)$ for the function spectrum.

Write $\lambda_N(k)$ for the real plane of the complex character
of weight $k$, with ordered basis $(1,i)$. For $C_{p^n}$, put
\[
\lambda_j=\lambda_{p^n}(p^j)\quad(0\le j<n).
\]
At $p=2$, $\sigma_{2^m}$ is the sign representation of $C_{2^m}$
for $m\ge1$, modeled on $\mathbb R$ with basis $1$.
The identification $\lambda_{2^n}(2^{n-1})=2\sigma_{2^n}$ sends
the ordered basis $(1,i)$ to the two ordered real basis vectors.
For $n\ge1$, the \emph{selected degrees} of \cite{ZhangHZ} form
\[
RO_{\mathrm{sel}}(C_{p^n})=
\begin{cases}
\mathbb Z\{1,\lambda_0,\ldots,\lambda_{n-1}\},&p\ne2,\\
\mathbb Z\{1,\sigma_{2^n},\lambda_0,\ldots,\lambda_{n-2}\},&p=2,
\end{cases}
\qquad RO_{\mathrm{sel}}(e)=\mathbb Z.
\]
\begin{definition}\label{def:degree-coordinates}
For $V\in RO(C_{p^n})$, its degree coordinates are
$a_i=\frac12\dim_{\mathbb R}V^{C_{p^i}}$ $(0\le i\le n)$.
For natural degrees, $\mathbf a=(a_0,\ldots,a_n)$ belongs to
\[
D_n^{(0)}=\mathbb Z^{n+1},\qquad
D_n^{(1)}=(\mathbb Z+\tfrac12)^n\times\mathbb Z\ (n\ge1),\qquad
D_n=\begin{cases}
D_n^{(0)},&p\ne2\text{ or }n=0,\\
D_n^{(0)}\sqcup D_n^{(1)},&p=2,\ n\ge1.
\end{cases}
\]
Set $D_0^{(1)}=\varnothing$ and $\epsilon(\mathbf a)=\epsilon$ for $\mathbf a\in D_n^{(\epsilon)}$.
Define
\begin{equation*}
\rho_p(\mathbf a)=\begin{cases}
2a_n+\displaystyle\sum_{j=0}^{n-1}(a_j-a_{j+1})\lambda_j,
   &p\ne2,\ n\ge1,\\
2a_n+2(a_{n-1}-a_n)\sigma_{2^n}
 +\displaystyle\sum_{j=0}^{n-2}(a_j-a_{j+1})\lambda_j,
   &p=2,\ n\ge1,\\
2a_0,&n=0.
\end{cases}
\end{equation*}
\end{definition}
For $L=C_{p^m}$, write $\rho_L$ for $\rho_p$ formed over $L$.
For $p=2$, $n\ge1$ and a natural degree $U$,
$2(a_{n-1}-a_n)\equiv\epsilon(\mathbf a)\pmod2$.
For $\epsilon(\mathbf a)=0$, we lift $U$ to a virtual complex
$\mathbb T$-representation whose restriction to $C_{p^n}$ has underlying
real degree $U$. For selected degrees, we use weights $p^j$ and restrict
the same $\mathbb T$-representation to each subgroup.
The Euler class is $a_W\in\pi_{-W}^G\mathbb S$. For $H\le G$, set
$a_H^G=a_{\mathbb R[G/H]-1}$, and use the same notation for its
images in $X_R$.

We put $x^+=\max(x,0)$. A line over $A$ is an invertible $A$-module.
For a line $L$ and an integer $a$, set
\[
L^a=\begin{cases}L^{\otimes a},&a\ge0,\\
(L^\vee)^{\otimes(-a)},&a<0.
\end{cases}
\]
Frobenius is $\varphi$, with $\varphi^*L=A\otimes_{A,\varphi}L$.
A \emph{graded unit} of degree $W$ has a multiplicative inverse of degree $-W$.

The rest of this paper is organized as follows.
Section~\ref{sec:compression} applies Tate compression to $X_R$.
For nonzero $p$-torsionfree perfectoid rings, Section~\ref{sec:uniform}
computes the graded Tambara functor and cyclotomic restriction.

Appendix~\ref{app:spectral-foundations} computes the local homotopy
groups and comparison maps used in Section~\ref{sec:uniform}.
Appendix~\ref{app:operation-foundations} gives the remaining
proofs for norms.
We treat perfect $\F_p$-algebras in Appendix~\ref{app:characteristic-p},
using the calculations for $\HZ$ in \cite{ZhangHZ}.

\section{Tate compression}\label{sec:compression}

In this section, we apply the Tate compression theorem to $X_R$. 

For $G=C_{p^n}$, set
$\mathcal F_i=\{H:H\le C_{p^i}\}$ and
$\mathcal F_{-1}=\varnothing$. For a genuine $G$-spectrum $X$, put
\[
X^{\Phi_i}=\widetilde E\mathcal F_{i-1}\wedge X,\qquad
X^{h_i}=F((E\mathcal F_i)_+,X^{\Phi_i}),\qquad
X^{t_i}=\widetilde E\mathcal F_i\wedge X^{h_i}.
\]
We use $X^{\Phi_i}$ for $0\le i\le n$ and $X^{h_i},X^{t_i}$ for $0\le i<n$.
Write $\Phi^{C_{p^i}}X$ for geometric fixed points and
$\theta_i:X^{\Phi_{i+1}}\to X^{t_i}$ for the map from isotropy separation~\cite[Section~17]{GreenleesMay}~\cite[Section~2.2]{AMR}.
A spectrum $Y$ is $\mathcal F_i$-complete if
$Y\to F((E\mathcal F_i)_+,Y)$ is an equivalence.

The following form of Tate compression applies to every representation
suspension of $X_R$.

\begin{proposition}\label{prop:compression}
Suppose each $\Phi^{C_{p^i}}X$ has bounded-below underlying spectrum.
For $0\le i<n-1$, $X^{t_i}$ is $\mathcal F_{i+1}$-complete, and
the map $\theta_i$ factors through
$X^{h_{i+1}}$ with contractible space of factorizations.
For $n\ge1$,
\begin{equation}\label{eq:tate-compression}
X\simeq\holim\bigl(
X^{h_0}\to X^{t_0}\leftarrow X^{h_1}
\to\cdots\to X^{t_{n-1}}\leftarrow X^{\Phi_n}\bigr).
\end{equation}
These statements apply to $X=S^V\wedge X_R$ for every $V\in RO(G)$.
\end{proposition}

\begin{proof}
For $r\ge i+1$, the completion map on $C_{p^r}$-fixed points is
\[
(\Phi^{C_{p^i}}X)^{tC_{p^{r-i}}}
\longrightarrow
\bigl((\Phi^{C_{p^i}}X)^{tC_p}\bigr)^{hC_{p^{r-i-1}}},
\]
an equivalence by \cite[Lemma~II.4.1]{NS}.
For $H\in\mathcal F_i$, the $H$-fixed points of the source and target
of the completion map are contractible. Hence
\[
X^{t_i}\xrightarrow{\simeq}
F((E\mathcal F_{i+1})_+,X^{t_i}).
\]
By the universal property of completion,
\[
\Map(X^{h_{i+1}},X^{t_i})\xrightarrow{\simeq}
\Map(X^{\Phi_{i+1}},X^{t_i}).
\]
Equation~\eqref{eq:tate-compression} is
\cite[Theorem~4.7]{ZhangHZ}.
On $G$-fixed points, this gives the iterated pullback of
\cite[Corollary~II.4.7 and Remark~II.4.8]{NS}.
Finally,
\[
\Phi^{C_{p^i}}(S^V\wedge X_R)
\simeq S^{V^{C_{p^i}}}\wedge\Phi^{C_{p^i}}X_R,
\qquad
\pi_k\Phi^{C_{p^i}}(S^V\wedge X_R)=0
\quad(k<\dim_{\mathbb R}V^{C_{p^i}}),
\]
since the geometric fixed points of $X_R$ are connective.
\end{proof}

The vertices of \eqref{eq:tate-compression} are $X^{h_i},X^{\Phi_n}$,
and its edges are $X^{t_i}$. For $n\ge1$, write
$\ell_i:X^{h_i}\to X^{t_i}$ ($0\le i<n$),
$\rho_i:X^{h_{i+1}}\to X^{t_i}$ ($0\le i<n-1$), and
$\rho_{n-1}=\theta_{n-1}:X^{\Phi_n}\to X^{t_{n-1}}$
for the maps in \eqref{eq:tate-compression}. The \emph{incidence map} is
\[
\Delta:\left(\prod_{i<n}X^{h_i}\right)\times X^{\Phi_n}
\longrightarrow\prod_{i<n}X^{t_i},\qquad
(\Delta x)_i=\ell_i(x_i)-\rho_i(x_{i+1}).
\]
Thus $X\simeq\operatorname{fib}(\Delta)$. For $n=0$, the diagram
is $X^{\Phi_0}=X$.

For $V\in RO(G)$, regard $V^{C_{p^i}}$ as a representation
of $C_{p^{n-i}}$ via the $p^i$-power identification. After taking
$G$-fixed points in \eqref{eq:tate-compression} for $S^{-V}\wedge X_R$,
the vertices and edges at $i<n$ are, respectively,
\[
(S^{-V^{C_{p^i}}}\wedge X_R)^{hC_{p^{n-i}}},\qquad
(S^{-V^{C_{p^i}}}\wedge X_R)^{tC_{p^{n-i}}}\quad(0\le i<n),
\]
and the terminal vertex is $S^{-V^{C_{p^n}}}\wedge X_R$. Each $\rho_i$ is induced by
\[
S^{-V^{C_{p^{i+1}}}}\wedge X_R
\simeq\Phi^{C_p}(S^{-V^{C_{p^i}}}\wedge X_R)
\longrightarrow(S^{-V^{C_{p^i}}}\wedge X_R)^{tC_p}.
\]
For $V=0$, this map is cyclotomic Frobenius.
For general $V$, its expression in terms of cyclotomic Frobenius
and the Euler map is \eqref{eq:app-moving-arrow}.

\section{The \texorpdfstring{$RO(C_{p^n})$}{RO(Cpn)}-graded Tambara functor}
\label{sec:uniform}
Let $R$ be a nonzero $p$-torsionfree perfectoid ring. We compute the
$RO(C_{p^n})$-graded Tambara functor of $X_R$.
The homotopy groups of the local terms are computed mainly
by combining earlier results. We give these computations in
Appendix~\ref{app:spectral-foundations}.
We first compute norms in natural degrees.
We then compute norms on the homotopy groups of $X_R^{t_i}$.
Finally, we use the boundary maps of \eqref{eq:tate-compression}
to compute norms in shifted degrees.

\subsection{The Breuil--Kisin twist}\label{subsec:period-line}
We need the following
Breuil--Kisin twist to describe the homotopy groups.

\begin{definition}\label{def:period-line}
Let $R$ be a nonzero $p$-torsionfree perfectoid ring with perfect prism
$(A,I_0,\varphi)$. Thus $A=W(R^\flat)$, where
$R^\flat=\varprojlim_{x\mapsto x^p}R/pR$, $I_0=\ker\theta$ for the
Fontaine map $\theta:A\to R$, and $\varphi$ is Witt Frobenius.
Put $I_j=\varphi^j(I_0)$. For $0\le i\le m$, write
$J_i^m=\prod_{i<j\le m}I_j$, with $J_i^i=A$. Define
\begin{equation*}
B=\operatorname{BK}(A,I_1)
   =\varprojlim_{r\ge1}\left(J_0^r/(J_0^r)^2,\,\tau_r\right).
\end{equation*}
The transition map $\tau_r:J_0^{r+1}/(J_0^{r+1})^2\to J_0^r/(J_0^r)^2$
is uniquely determined by
\[
p\tau_r([x])=[x]\qquad(x\in J_0^{r+1}).
\]
Existence and uniqueness follow from
\cite[Corollary~2.2.10]{BhattLurie}.
\end{definition}

By \cite[Constructions 2.2.11, 2.2.14 and
Proposition 2.2.12]{BhattLurie}, $B$ is an invertible $A$-module and
\begin{equation}\label{eq:algebraic-period-phi}
 \Phi_B:\varphi^*B\xrightarrow{\sim}I_1^{-1}B,
 \qquad \varphi(J_0^r)=I_1^{-1}J_0^{r+1}.
\end{equation}
Write
\[
\operatorname{can}:
\operatorname{TC}^-(R;\Z_p)=X_R^{h\mathbb T}
\longrightarrow X_R^{t\mathbb T}=\operatorname{TP}(R;\Z_p).
\]
Use the complex orientation under
$\mathbb T/C_{p^r}\xrightarrow{\sim}\mathbb T$,
$[z]\mapsto z^{p^r}$.
By \cite[Proposition~6.5 and Remark~6.6]{BMS},
\[
\begin{tikzcd}[column sep=large,row sep=small]
\pi_2\operatorname{TC}^-(R;\Z_p)
  \arrow[r,"\pi_2(\operatorname{can})"] \arrow[d,"\cong"']
& \pi_2\operatorname{TP}(R;\Z_p) \arrow[d,"\cong"] \\
I_0B \arrow[r,hook] & B.
\end{tikzcd}
\]
For $r\ge1$, the cyclotomic map
$X_R^{hC_{p^{r-1}}}\to X_R^{tC_{p^r}}$
is $\mathbb T$-equivariant.
In the homotopy fixed-point spectral sequences for the residual circle,
use its complex Euler class to identify $\pi_2$ with the first
quotient. Naturality \cite[Remark~6.6]{BMS} gives
\[
\begin{tikzcd}[column sep=large,row sep=large]
\pi_2 X_R^{hC_{p^{r-1}}}
  \arrow[r] \arrow[d,"\cong"']
& \pi_2 X_R^{tC_{p^r}} \arrow[d,"\cong"] \\
\dfrac{I_0J_0^{r-1}}{(I_0J_0^{r-1})^2}
  \arrow[r,"\varphi"']
& \dfrac{J_0^r}{(J_0^r)^2}.
\end{tikzcd}
\]
The inverse limit and \cite[Construction~2.2.14]{BhattLurie}
for $(A,I_1)$ give the degree-two linearization of cyclotomic Frobenius:
\[
\varphi^*(I_0B)=I_1\varphi^*B
 \xrightarrow{\ I_1\Phi_B\ }B.
\]
We use the canonical isomorphism
\begin{equation}\label{eq:period-cotangent-reduction}
B/I_0B\xrightarrow{\sim}R,\qquad
I_0B/I_0^2B\xrightarrow{\sim}I_0/I_0^2
\end{equation}
of \cite[Proposition~6.5 and the preceding discussion]{BMS}. By \cite[Theorem~6.1]{BMS},
$\pi_{2a}X_R\cong(I_0/I_0^2)^a$ for $a\ge0$. Projection to $r=1$ and iteration of \eqref{eq:algebraic-period-phi} give
\begin{equation}\label{eq:period-repair-data}
B/I_1B\simeq I_1/I_1^2,\qquad
 \Sigma_r:\varphi^{r*}B\xrightarrow{\sim}(J_0^r)^{-1}B\quad(r\ge0).
\end{equation}
When $I_0$ is principal, we can use the first isomorphism
in \eqref{eq:period-repair-data} to choose a generator of $B$.

\begin{proposition}\label{lem:period-repair-frame}
If $I_0=(\xi)$, every lift $b\in B$ of
$\varphi(\xi)\bmod I_1^2$ generates $B$ freely over $A$.
\end{proposition}
\begin{proof}
The ring $A$ is $(p,I_1)$-complete, so $I_1\subseteq\operatorname{Jac}(A)$.
The map $A\to B$, $1\mapsto b$, is an isomorphism modulo $I_1$.
By Nakayama's lemma it is surjective. Since $B$ is invertible, it is
an isomorphism.
\end{proof}

\subsection{Homotopy groups in representation degrees}
We first reduce to $RO_{\mathrm{sel}}(C_{p^n})$
using maps of representation spheres.

\begin{definition}\label{def:character-translation}
For $0\le j<n$, $r>0$, and $(r,p)=1$, the map
\[
f_{r,j}:S^{\lambda_{p^n}(p^j)}
\longrightarrow S^{\lambda_{p^n}(rp^j)},
\qquad f_{r,j}(z)=z^r,
\]
is a $p$-local equivalence. Its map on $C_{p^l}$-fixed spheres has
degree $r$ for $l\le j$ and degree $1$ for $l>j$. Define the unit
\[
\widetilde v_{r,j}=r f_{r,j}^{-1}
\in
\pi_{\lambda_{p^n}(rp^j)-\lambda_{p^n}(p^j)}
^{C_{p^n}}\mathbb S_{(p)},
\]
with fixed-point degrees
\[
\deg(\widetilde v_{r,j}^{C_{p^l}})
=
\begin{cases}
1,&l\le j,\\
r,&l>j.
\end{cases}
\]
For
\[
D=\sum_{r,j}m_{r,j}
\bigl(\lambda_{p^n}(rp^j)-\lambda_{p^n}(p^j)\bigr),
\qquad m_{r,j}\in\mathbb Z,
\]
define
\[
\widetilde v_D
=\prod_{r,j}\widetilde v_{r,j}^{m_{r,j}}.
\]
The unit $\widetilde v_D$ depends on the weights $r$. Under restriction we use the same maps $f_{r,j}$.
Let $v_{r,j},v_D$ be the images of $\widetilde v_{r,j},\widetilde v_D$ in $X_R$.
For $V=V_0+D$ with $V_0$ selected, multiplication
by $v_D$ induces an isomorphism
\[
\pi_{V_0}^{C_{p^n}}X_R
\xrightarrow{\ \cong\ }
\pi_V^{C_{p^n}}X_R.
\]
\end{definition}

For $V\in RO(C_{p^n})$, obtain $V_0$ by replacing
$\lambda_{p^n}(rp^j)$ by $\lambda_j$ in its irreducible
expansion, where $0<rp^j<p^n/2$ and $(r,p)=1$,
and leaving the one-dimensional summands unchanged.
Use the expansion of $D=V-V_0$ in
Definition~\ref{def:character-translation}.
The local homotopy groups in degrees $V_0$ and $V$ are
identified by multiplication by $v_D$,
which commutes with $\ell_i$ and $\rho_i$.

We have the following definition.

\begin{definition}\label{def:natural-degree}\label{def:operation-complexes}
Let $R$ be a nonzero $p$-torsionfree perfectoid ring.
Fix $n\ge0$ and $V\in\RO(C_{p^n})$, and put
\[
U=\begin{cases}V,&\dim V^{C_{p^n}}\equiv0\pmod2,\\
V+1,&\dim V^{C_{p^n}}\equiv1\pmod2,
   \end{cases}
\qquad a_i=\frac12\dim_{\mathbb R}U^{C_{p^i}}\quad(0\le i\le n).
\]
For $p$-torsionfree $R$, define $A_i(U)$ and $T_i(U)$ for $0\le i<n$
as follows. Products of invertible ideals and their
negative powers are taken in the localization of $A$ obtained by
inverting their generators.

If $\epsilon(\mathbf a)=0$, set
\begin{align}
E_i(U)&=(\varphi^{n*}B)^{a_n}\prod_{i<j\le n}I_j^{a_j}
\quad(0\le i\le n),\notag\\
A_i(U)&=\frac{I_i^{a_i^+}E_i(U)}
                {J_i^n I_i^{(a_i+1)^+}E_i(U)},
 &T_i(U)&=E_i(U)/J_i^n E_i(U).\label{eq:uniform-complex}
\end{align}

If $\epsilon(\mathbf a)=1$, set
\begin{align}
\mathcal S_i(U)&=(\varphi^{n*}B)^{a_n}I_n^{a_n}\prod_{i<j<n}I_j^{a_j+1/2},\notag\\
A_i(U)&=I_i^{(a_i+1/2)^+}\mathcal S_i(U)/I_nI_i^{(a_i+1/2)^+}\mathcal S_i(U),
 &T_i(U)&=\mathcal S_i(U)/I_n\mathcal S_i(U).\label{eq:uniform-sign}
\end{align}
For $p$-torsionfree $R$, set
\[
A_n(U)=\begin{cases}
\displaystyle
\frac{I_n^{a_n}(\varphi^{n*}B)^{a_n}}
     {I_n^{a_n+1}(\varphi^{n*}B)^{a_n}},&a_n\ge0,\\[6pt]
0,&a_n<0.
\end{cases}
\]
For $0\le i<n$, define $\ell_i:A_i(U)\to T_i(U)$ and
$\rho_i:A_{i+1}(U)\to T_i(U)$ as follows.
For $p$-torsionfree $R$, both maps send $[x]$ to $[x]$, with
representatives in $E_i(U)$ if $\epsilon(\mathbf a)=0$ and in
$\mathcal S_i(U)$ if $\epsilon(\mathbf a)=1$.
For $0\le i<n-1$, use the identifications
\[
E_i(U)=I_{i+1}^{a_{i+1}}E_{i+1}(U),
\qquad
\mathcal S_i(U)=I_{i+1}^{a_{i+1}+1/2}\mathcal S_{i+1}(U).
\]
For $n\ge1$ and $a_n\ge0$, we have $A_n(U)=T_{n-1}(U)$ and
$\rho_{n-1}=\mathrm{id}$; for $a_n<0$, the source is zero.

For $G=C_{p^n}$, we use the incidence map to define the following complex.
\[
\begin{gathered}
C_G^\bullet(U)=
\left[\bigoplus_{i=0}^n A_i(U)\xrightarrow{\Delta_U}
\bigoplus_{i=0}^{n-1}T_i(U)\right],\\
(\Delta_Uq)_i=\ell_i(q_i)-\rho_i(q_{i+1}),
\end{gathered}
\]
with cochain degrees $0,1$.
\end{definition}

Sulyma computes finite-level $\TR$ groups in a
contiguous range of representation degrees~\cite[Theorem~4.12]{SulymaNew}.
For every natural degree $U$, the kernel and cokernel of $\Delta_U$ give
the homotopy groups in degrees $U$ and $U-1$, respectively.

\begin{proposition}\label{thm:uniform}
For $R$, $n$ and $U$ as in Definition~\ref{def:natural-degree},
\begin{equation}\label{eq:uniform-groups}
 \pi_U^{C_{p^n}}X_R=\ker\Delta_U,
 \qquad \pi_{U-1}^{C_{p^n}}X_R=\coker\Delta_U.
\end{equation}
\end{proposition}

\begin{proof}
For $U\in RO_{\mathrm{sel}}(C_{p^n})$, the identifications
\[
A_i(U)=\pi_U^{C_{p^n}}X_R^{h_i},\qquad
T_i(U)=\pi_U^{C_{p^n}}X_R^{t_i},\qquad
A_n(U)=\pi_U^{C_{p^n}}X_R^{\Phi_n}
\]
with the first two taken for $0\le i<n$, follow from
Proposition~\ref{prop:app-finite-comparison},
\eqref{eq:app-frob-numerator}--\eqref{eq:app-frob-denominator}, and
\eqref{eq:app-sign-local-kernels}.
Multiplication by $v_D$ gives the same identifications
for all natural degrees. The same computations give
\[
\pi_{U\pm1}^{C_{p^n}}X_R^{h_i}=0,\qquad
\pi_{U\pm1}^{C_{p^n}}X_R^{t_i}=0,\qquad
\pi_{U\pm1}^{C_{p^n}}X_R^{\Phi_n}=0.
\]
The long exact sequence of the fiber of the incidence map
in Section~\ref{sec:compression} gives \eqref{eq:uniform-groups}.
\end{proof}

\subsection{Products, restrictions, transfers, and Weyl actions}
\label{subsec:products-maps}

We now give the formulas for products, restrictions, transfers, and Weyl actions.

When comparing $G=C_{p^n}$ with $H=C_{p^{n-1}}$, insert a zero
summand at $A_n$ and $T_{n-1}$ in $C_H^\bullet$.

For natural degrees $U,W$ with degree coordinates
$\mathbf a,\mathbf b$ satisfying $\epsilon(\mathbf a)=\epsilon(\mathbf b)=0$,
the products at $A_i$ and $T_i$ for $i<n$ are induced by
$E_i(U)\otimes_A E_i(W)\simeq E_i(U+W)$.

At $p=2$, when $\epsilon(\mathbf a)=1$, we use the sign Euler class to
express restriction and transfer as connecting homomorphisms.

\begin{proposition}\label{prop:sign-exact-complexes}
Let $G=C_{2^n}$, $H=C_{2^{n-1}}$, and $n\ge1$.
If $U$ is natural with $\epsilon(\mathbf a)=1$,
there are degreewise exact sequences
\begin{equation}\label{eq:sign-exact-complexes}
\begin{aligned}
0\longrightarrow C_H^\bullet(\Res_H^G(U+\sigma_{2^n}))
&\xrightarrow{\tr}C_G^\bullet(U+\sigma_{2^n})
\xrightarrow{a_{\sigma_{2^n}}}C_G^\bullet(U)\longrightarrow0,\\
0\longrightarrow C_G^\bullet(U)
&\xrightarrow{a_{\sigma_{2^n}}}C_G^\bullet(U-\sigma_{2^n})
\xrightarrow{\res}C_H^\bullet(\Res_H^G(U-1))\longrightarrow0.
\end{aligned}
\end{equation}
Both maps $a_{\sigma_{2^n}}$ send $[q]$ to $[q]$.
The connecting homomorphisms on cohomology are
\[
\begin{aligned}
\res_H^G&:H^0(C_G^\bullet(U))
\longrightarrow H^1(C_H^\bullet(\Res_H^G(U+\sigma_{2^n}))),\\
\tr_H^G&:H^0(C_H^\bullet(\Res_H^G(U-1)))
\longrightarrow H^1(C_G^\bullet(U)).
\end{aligned}
\]
\end{proposition}

\begin{proof}
Firstly we apply $(G/H)_+\to S^0\to S^{\sigma_{2^n}}$ and its dual to
\eqref{eq:tate-compression}. By Shapiro's equivalence and the vanishing
in shifted degrees, we have  \eqref{eq:sign-exact-complexes}.

Write $\sigma=\sigma_{2^n}$. For $i<n$, the coefficient maps are
\[
A_i(U+\sigma_{2^n})\twoheadrightarrow A_i(U)
\hookrightarrow A_i(U-\sigma_{2^n}),
\]
and likewise for $T_i$.
Indeed, $E_i(U+\sigma)=\mathcal S_i(U)=J_i^{n-1}E_i(U-\sigma)$.
The second map is the kernel inclusion in \eqref{eq:app-sign-local-kernels}.
Their composite is $a_{2\sigma}=a_{\lambda_{n-1}}$, which acts by inclusion
on the circle representatives \cite[Theorem~1.9 and Example~3.2]{SulymaNew}.
Hence the first map sends $[q]$ to $[q]$ as well. At $A_n$, $\sigma^G=0$.

By naturality,  the connecting homomorphisms are $\res_H^G$
and $\tr_H^G$.
\end{proof}

For $U=\rho_p(\mathbf a)$ with $\mathbf a\in D_n^{(0)}$, we use the same virtual complex
$\mathbb T$-representation lifting $U$ to compute the terms at
$G$ and at $H$.
For $V=V_0+D$, we can transport the sphere maps by
\begin{equation}\label{eq:character-res-tr}
\begin{aligned}
\res_H^G(v_Dx_0)
 &=\res_H^G(v_D)\res_H^G(x_0),\\
\tr_H^G(y)
 &=v_D\,\tr_H^G\!\left(\res_H^G(v_D)^{-1}y\right).
\end{aligned}
\end{equation}
Here $x_0\in\pi_{V_0}^G X_R$ and
$y\in\pi_{\Res_H^G V}^H X_R$.
The above identities follow from naturality and the projection formula.
For $p$-torsionfree $R$, \cite[Proposition~2.2.12]{BhattLurie}
for $(A,I_1)$ and \eqref{eq:algebraic-period-phi} give
\[
 B/(I_0J_0^r)B\simeq J_0^r/(I_0J_0^r)J_0^r,\qquad
 B/J_0^rB\simeq J_0^r/(J_0^r)^2\quad(r\ge0).
\]
The first follows from $B/J_0^{r+1}B\simeq J_0^{r+1}/(J_0^{r+1})^2$
by pulling back along $\varphi^{-1}$ and tensoring with $I_0^{-1}$.
For $U=\rho_p(\mathbf a)$ with $\mathbf a\in D_n^{(0)}$, the Frobenius isomorphism gives
\[
E_i(U)\cong
(\varphi^{(n-1)*}B)^{a_n-a_{n-1}}\otimes_A E_i(\Res_H^G U),
\qquad 0\le i<n,
\]
where the two lines are formed at $G$ and $H$, respectively.
The quotient isomorphisms and $\varphi^{i*}\Sigma_{n-1-i}$ identify
$\varphi^{(n-1)*}B$ with $A$ modulo $I_iJ_i^{n-1}$.
Thus, using the same circle lift, the terms at $H$ are
\begin{equation}\label{eq:subgroup-coefficients}
\begin{aligned}
A_i(\Res_H^G U)&\cong
\frac{I_i^{a_i^+}E_i(U)}
     {J_i^{n-1}I_i^{(a_i+1)^+}E_i(U)}
&& (0\le i<n),\\
T_i(\Res_H^G U)&\cong E_i(U)/J_i^{n-1}E_i(U)
&& (0\le i<n-1).
\end{aligned}
\end{equation}

We now define the transfer coefficient.

\begin{definition}\label{def:transfer-coefficient}
For $m\ge0$, let
$V:W_{m+1}(R)\to W_{m+2}(R)$
be the Witt vector Verschiebung.
For $p$-torsionfree $R$, the maps
$\widetilde\theta_r:A\to W_r(R)$
of \cite[Lemma~3.12]{BMS18} give
\[
\widetilde\theta_{m+1}\varphi:
A/(I_0J_0^m)\xrightarrow{\sim}W_{m+1}(R),
\qquad
v_m=(\widetilde\theta_{m+2}\varphi)^{-1}(V(1))
\in A/(I_0J_0^{m+1}).
\]
\end{definition}

We can now give the formulas.

\begin{proposition}\label{prop:products-maps}\label{lem:path-products}
Let $R$ be as in Definition~\ref{def:natural-degree} and
let $G=C_{p^n}$.

\begin{enumerate}[label=\textup{(\arabic*)},leftmargin=*]
\item For natural degrees $U,W$ with degree coordinates
$\mathbf a,\mathbf b$, respectively, let
\[
x\in\ker\Delta_U,\quad y\in\ker\Delta_W,\quad
[z]\in\coker\Delta_U,\quad [w]\in\coker\Delta_W.
\]
The unit in degree zero is $(1,\ldots,1)$.
Under \eqref{eq:uniform-groups}, the products are
\begin{equation}\label{eq:path-cup}
\begin{aligned}
(xy)_j&=x_jy_j &&(0\le j\le n),\\
x[w]&=[(\ell_i(x_i)w_i)_{i<n}],\\
[z][w]&=0\in\pi_{U+W-2}^G X_R.
\end{aligned}
\end{equation}

\item The restriction and transfer maps on kernels and cokernels are as follows.
Let $n\ge1$ and $H=C_{p^{n-1}}$. For a natural degree $U$
with $\epsilon(\mathbf a)=0$, the maps
\[
\begin{aligned}
\res_H^G&:C_G^\bullet(U)\longrightarrow C_H^\bullet(\Res_H^G U),\\
\tr_H^G&:C_H^\bullet(\Res_H^G U)\longrightarrow C_G^\bullet(U)
\end{aligned}
\]
are given, for $p$-torsionfree $R$ and selected $U$, in the
coordinates of \eqref{eq:subgroup-coefficients}, by
\begin{equation}\label{eq:adjacent-res-tr}
\begin{aligned}
\res_H^G(x_0,\ldots,x_n)&=( x_0,\ldots, x_{n-1}),\\
\tr_H^G(y_0,\ldots,y_{n-1})&=(v_{n-1}y_0,\ldots,v_{n-1}y_{n-1},0),\\
\res_H^G(z_0,\ldots,z_{n-1})&=( z_0,\ldots, z_{n-2}),\\
\tr_H^G(w_0,\ldots,w_{n-2})&=(v_{n-1}w_0,\ldots,v_{n-1}w_{n-2},0).
\end{aligned}
\end{equation}
For arbitrary $U$, use \eqref{eq:character-res-tr}.

\item Let $p=2$, $n\ge1$, and $H=C_{2^{n-1}}$.
Suppose $U$ is natural with
$\epsilon(\mathbf a)=1$. For
\[
x\in\ker\Delta_U,\qquad
y\in\ker\Delta_{\Res_H^G(U-1)},
\]
choose lifts
\[
b\in C_G^0(U+\sigma_{2^n}),\quad a_{\sigma_{2^n}} b=x,
\qquad
c\in C_G^0(U-\sigma_{2^n}),\quad\res_H^G(c)=y.
\]
Then
\begin{equation}\label{eq:sign-res-tr}
\begin{aligned}
\res_H^G(x)&=[\tr^{-1}\Delta_{U+\sigma_{2^n}}b]
\in\coker\Delta_{\Res_H^G(U+\sigma_{2^n})},\\
\tr_H^G(y)&=[a_{\sigma_{2^n}}^{-1}\Delta_{U-\sigma_{2^n}}c]
\in\coker\Delta_U.
\end{aligned}
\end{equation}
Moreover,
\[
\tr_H^G\big|_{\pi_{\Res_H^G U}^H X_R}=0,
\qquad
\res_H^G\big|_{\pi_{U-1}^G X_R}=0.
\]
\end{enumerate}
\end{proposition}

\begin{proof}
\textup{(1)} The maps in diagram of \eqref{eq:tate-compression}
are maps of commutative $X_R$-algebras in $\Sp^G$, so the
associated fiber sequence is a sequence of $X_R$-modules.
By the proof of \cite[Theorem~4.30, equations~(4.31)--(4.32)]{ZhangHZ} and \eqref{eq:uniform-groups}, we have the three product
formulas in \eqref{eq:path-cup}.

For $p$-torsionfree $R$, tensoring the ideals in
\eqref{eq:uniform-complex} and \eqref{eq:uniform-sign}
gives the products formula.
For $\epsilon(\mathbf a)=\epsilon(\mathbf b)=1$ and $0\le i<n$,
multiplication is induced by the inclusion
\[
\mathcal S_i(U)\mathcal S_i(W)
=J_i^{n-1}E_i(U+W)\subseteq E_i(U+W).
\]
The unit has coordinate $1$ at every vertex.

\textup{(2)} Restriction and transfer commute with the incidence map
and its connecting homomorphism. Thus they induce the maps on
kernels and cokernels, as in the proof of
\cite[Theorem~6.13]{ZhangHZ}.
For $p$-torsionfree $R$,
\cite[Proposition~6.2 and Remark~6.6]{BMS} and
\cite[Theorem~3.3]{HM97} give
\[
\tr_H^G(1)
=(\widetilde\theta_{n+1}\varphi)^{-1}(V(1))
=v_{n-1}.
\]
Restriction is the quotient in \eqref{eq:subgroup-coefficients} at each vertex and edge.
Transfer follows from the projection formula, proving \eqref{eq:adjacent-res-tr}.

\textup{(3)} The proofs of
\cite[Propositions~6.17--6.18]{ZhangHZ}, applied to
\eqref{eq:sign-exact-complexes}, give
\eqref{eq:sign-res-tr} up to sign.
Since $X_R$ has a circle action and $\epsilon(\mathbf a)=1$,
\[
\begin{aligned}
\res_H^G(x)
 &=c_{\gamma_{2^n}}\res_H^G(x)=-\res_H^G(x),\\
\tr_H^G(y)
 &=\tr_H^G(c_{\gamma_{2^n}}y)=-\tr_H^G(y).
\end{aligned}
\]
The proof of \cite[Proposition~6.20]{ZhangHZ} gives
$\tr_H^G|_{\pi_{\Res_H^G U}^H X_R}=0$ and
$\res_H^G|_{\pi_{U-1}^G X_R}=0$:
the map
\[
X_R\longrightarrow
\left(\prod_{i<n}X_R^{h_i}\right)\times X_R^{\Phi_n}
\]
induces zero on the source groups and an injection on the target groups
by Proposition~\ref{thm:uniform}.
\end{proof}

It remains to determine how a generator of $G$ acts on
$\pi_{\Res_L^G V}^L X_R$ for $L<G$.

\begin{proposition}
\label{prop:weyl-action}
Let $G=C_{p^n}$, $V\in RO(G)$, and let $L<G$ be any proper
subgroup. Let $\gamma$ be the generator of $G$.
For odd $p$, the Weyl action on
$\pi_{\Res_L^G V}^L X_R$ is trivial.
For $p=2$, write
\[
V=m+s\sigma_{2^n}+\sum_\lambda s_\lambda\lambda,
\]
where $\lambda$ ranges over the irreducible real rotation
representations. Then
\[
c_{\gamma^a}\big|_{\pi_{\Res_L^G V}^L X_R}
=(-1)^{as}\operatorname{id}
\qquad(a\in\mathbb Z).
\]
\end{proposition}

\begin{proof}
For $f\in[S^{\Res_L^G V},X_R]^L$,
\[
c_\gamma(f)
=\gamma_{X_R}\circ f\circ\gamma_{S^V}^{-1}.
\]
A path in $\mathbb T$ from $1$ to $\gamma$ gives
$\gamma_{X_R}\simeq_L\operatorname{id}$.
For $p=2$,
\[
L\leq\ker\sigma_{2^n},\qquad
\gamma_{S^{\sigma_{2^n}}}\simeq_L-\operatorname{id},
\qquad
\gamma_{S^\lambda}\simeq_L\operatorname{id},
\]
where the last homotopy is through rotations.
The sphere identities extend to virtual summands, so
\[
\gamma_{S^V}\simeq_L(-1)^s\operatorname{id},
\qquad
c_{\gamma^a}(f)=(-1)^{as}f.
\]
For odd $p$, every nontrivial irreducible real representation
is a rotation, so we have $c_\gamma=\operatorname{id}$.
\end{proof}

\subsection{Norms}\label{subsec:graded-norms}
We now compute the norms.

Let $H=C_{p^h}<G=C_{p^n}$. We first compute
\[
N_H^G:\pi_V^H X_R\longrightarrow\pi_{\Ind_H^G V}^G X_R,
\qquad V\in RO_{\mathrm{sel}}(H).
\]

\subsubsection{Fixed-point coordinates of induced representations}\label{subsubsec:selected-degree-norms}
For $V\in RO(H)$, write
\[
\begin{aligned}
2a_j&=\dim V^{C_{p^j}} &&(0\le j\le h),\\
2A_j&=\dim(\Ind_H^G V)^{C_{p^j}} &&(0\le j\le n).
\end{aligned}
\]

To determine the target degree of $N_H^G$, we compute the $A_j$
from the $a_j$.

\begin{proposition}\label{prop:induction-of-degrees}
The degree coordinates of $\Ind_H^G V$ satisfy
\begin{equation}\label{eq:period-repair-induced-degree}
A_j=\begin{cases}
p^{n-h}a_j,&0\le j\le h,\\
p^{n-j}a_h,&h<j\le n.
\end{cases}
\end{equation}
In particular, for $n=h+1$,
\begin{equation*}
I(\mathbf a)=(pa_0,\ldots,pa_h,a_h).
\end{equation*}
If $V$ is natural, then $\Ind_H^G V$ is natural and
$\epsilon(\mathbf A)=0$.
\end{proposition}

\begin{proof}
For $K=C_{p^j}$, the fixed subspace of the induced representation is
\[
(\Ind_H^G V)^K\cong
\bigoplus_{K\backslash G/H}V^{K\cap H}.
\]
Here $K\cap H=C_{p^{\min(j,h)}}$ and
$|K\backslash G/H|=p^{n-\max(j,h)}$, which gives
\eqref{eq:period-repair-induced-degree}. The same formula holds
for virtual representations by additivity.
\end{proof}

\subsubsection{Norms in natural degrees}
Let $U\in RO_{\mathrm{sel}}(H)$ be natural.
We use geometric fixed points to determine $N_H^G$ on $\pi_U^H X_R$.

We first give the following definition.

\begin{definition}\label{def:geometric-fixed-point-maps}
For $U\in RO(C_{p^n})$, define
\[
\begin{aligned}
\gamma_j\colon
\pi_U^{C_{p^n}}X_R
&\xrightarrow{\Phi^{C_{p^j}}}
\pi_{U^{C_{p^j}}}^{C_{p^{n-j}}}
\Phi^{C_{p^j}}X_R\\
&\xrightarrow[\cong]{\mathrm{cyclotomic}}
\pi_{U^{C_{p^j}}}^{C_{p^{n-j}}}X_R
\xrightarrow{\operatorname{res}_e^{C_{p^{n-j}}}}
\pi_{U^{C_{p^j}}}^{e}X_R,
\qquad 0\le j\le n.
\end{aligned}
\]
For \(m\ge 1\), let
\[
\begin{gathered}
[m]\colon S^1\longrightarrow S^1,
\qquad z\longmapsto z^m,\\
\psi^m\colon X_R\longrightarrow [m]^*X_R,
\end{gathered}
\]
in $\Sp^{B\mathbb T}$, where $[m]^*$ pulls back the circle action along $[m]$ and
$\psi^m$ is the Adams operation induced by $[m]$
\cite[Definition~5.1 and Remark~5.2]{SulymaNew}.
\end{definition}

For $p$-torsionfree $R$, the Adams operations satisfy the following.
On $\pi_*X_R$, we have $\psi^{p^r}=(\psi^p)^r$.
By \cite[Remark~5.2 and Lemma~5.3]{SulymaNew},
\begin{equation}\label{eq:rcyc-crossing}
R_{\mathrm{cyc}}N_e^{C_p}=\psi^p,
\qquad
\left.\psi^{p^r}\right|_{\pi_{2d}X_R}
=p^{rd}\operatorname{id}
\quad(r,d\ge0).
\end{equation}
The first identity uses pullback along $\varrho_p$ and,
for a finite-dimensional real vector space $W$, the identification
\[
W\xrightarrow{\sim}(\Ind_e^{C_p}W)^{C_p},
\qquad w\longmapsto(w,\ldots,w).
\]
For negative summands, use the duals of the inverse sphere maps.

To recover a norm from its images under the maps $\gamma_j$, we need
to know when these maps detect a class.

\begin{lemma}\label{lem:geometric-ghost-injectivity}
Let $R$ be $p$-torsionfree and let $U\in RO(H)$ be natural, with
degree coordinates $\mathbf a=(a_0,\ldots,a_h)$. The maps
\[
\begin{aligned}
(\gamma_0,\ldots,\gamma_h):\pi_U^H X_R
&\longrightarrow\prod_{j=0}^h\pi_{U^{C_{p^j}}}^eX_R,
&&\epsilon(\mathbf a)=0,\\
\gamma_h:\pi_U^H X_R&\longrightarrow\pi_{U^H}^eX_R,
&&\epsilon(\mathbf a)=1,
\end{aligned}
\]
are injective. 
\end{lemma}
\begin{proof}
Multiplication by $v_D$ reduces the assertion to
$U\in RO_{\mathrm{sel}}(H)$.
By \cite[Lemma~2.2.8]{BhattLurie},
$\theta\varphi^{-i}(I_j)=pR$ for $j>i$.
For an invertible $A$-module $E$, a local generator of $I_j$
therefore acts injectively on $E/I_iE$, since $R$ is $p$-torsionfree.
The $I_i$-adic filtration gives injectivity on $E/I_i^bE$ for every
$b\ge0$, also for a local generator of $J_i^h$. Hence
\begin{equation}\label{eq:ghost-intersection}
\begin{aligned}
I_i^bE\cap I_jE&=I_i^bI_jE,\qquad 0\le i<j\le h,\\
I_i^bE\cap J_i^hE&=I_i^bJ_i^hE,\qquad 0\le i<h.
\end{aligned}
\end{equation}

Suppose $\epsilon(\mathbf a)=0$ and
$x=(x_0,\ldots,x_h)\in\ker\Delta_U$ satisfies $\gamma_jx=0$ for
$0\le j\le h$. Cyclotomic equivalence gives $x_h=0$.
If $x_{i+1}=0$, then
\[
\ell_i(x_i)=\rho_i(x_{i+1})=0.
\]
For $a_i\ge0$, \eqref{eq:ghost-line-conversion} and $\gamma_ix=0$
place a representative of $x_i$ in $I_i^{a_i+1}E_i(U)$.
Since $\ell_i(x_i)=0$, this representative also lies in $J_i^hE_i(U)$.
Equation~\eqref{eq:ghost-intersection} gives $x_i=0$.
For $a_i<0$, $\ell_i$ is an isomorphism.
Induction gives $x=0$.

Suppose $\epsilon(\mathbf a)=1$.
Equation~\eqref{eq:ghost-intersection} with $j=h$ implies that
the maps $\ell_i$ of \eqref{eq:uniform-sign} are injective for $0\le i<h$.
If $\gamma_hx=0$, then $x_h=0$ and
$\ell_i(x_i)=\rho_i(x_{i+1})$ gives $x_i=0$ for $i=h-1,\ldots,0$.

\end{proof}

Let $M$ be the ordered model of $\Ind_H^G U$ and
$\eta:S^M\xrightarrow{\sim}S^{\Ind_H^G U}$ its grading map
\cite[Definition~5.1 and Remark~5.2]{GradedTambara}.
For each $\chi:G\to\{\pm1\}$, by $\sqrt{[G:H]}$ times evaluation at $H$
we identify the $\chi$-isotypic subspace of $\Ind_H^G U$ with the ordered
$\chi|_H$-isotypic subspace of $U$. For $K=C_{p^j}$ let
\[
d_j:\bigoplus_{O\in K\backslash G/H} U^{K\cap H}
  \xrightarrow{\sim}(\Ind_H^G U)^K
\]
be the orbit diagonal. Order orbits by their least cosets
in the order $H,\gamma_{p^n}H,\ldots$, and set
$d_j(v)=|O|^{-1/2}v$ at the least coset of $O$.

\begin{proposition}\label{prop:vertex-norms}
Let $R$ be $p$-torsionfree and let $U\in RO_{\mathrm{sel}}(H)$ be natural.
Write $M=M_0+D$ as in Definition~\ref{def:character-translation}.
For $x\in\pi_U^H X_R$, the element
$y=v_D^{-1}\eta^*N_H^Gx\in\ker\Delta_{M_0}$ is given by
\begin{equation*}
\gamma_j(y)=\frac{\deg(d_j^{-1}\Phi^{C_{p^j}}\eta)}{\gamma_j(v_D)}
\begin{cases}
(\gamma_jx)^{p^{n-h}},&0\le j\le h,\\
\bigl(\psi^{p^{j-h}}\gamma_hx\bigr)^{p^{n-j}},&h<j\le n.
\end{cases}
\end{equation*}
\end{proposition}

\begin{proof}
In the orbit basis $d_j$, $\Phi^K$ gives $[G:HK]$ factors
\cite[Appendix~A]{HHR}, each equal to $\gamma_jx$ for $j\le h$
and $\psi^{p^{j-h}}\gamma_hx$ for $j>h$
by \eqref{eq:rcyc-crossing}. Pullback along $\Phi^K\eta$ contributes
$\deg(d_j^{-1}\Phi^K\eta)$. Uniqueness follows from
Proposition~\ref{prop:induction-of-degrees} and
Lemma~\ref{lem:geometric-ghost-injectivity}.
\end{proof}

In degree zero this recovers the Witt norm
\cite[Theorems~1.4 and~4.1]{WittNorm} and its prismatic form
\cite[Section~3.1]{SulymaPrisms}.

For products and composites, go from the induced grading map
$\eta'$ to $\eta$ acts by the sphere unit $(\eta')^{-1}\eta$.

\subsubsection{Norms on $X_R^{t_i}$}
For $0\le i<h<n$, Proposition~\ref{prop:boundary-indexed-multiplications} constructs
compatible $G/H$-indexed multiplications on the spectra in the
Tate compression diagram truncated at $X_R^{\Phi_{i+1}}$.
For each spectrum $Y$ in this diagram, write
\[
\mu_Y:N_H^G\res_H^G Y\longrightarrow Y.
\]
For each arrow $f:Y\to Z$, we have
\[
f\mu_Y=\mu_ZN_H^G(\res_H^Gf).
\]
For $U\in RO(H)$, write
\[
N_H^G:\pi_U^H X_R^{t_i}
\longrightarrow\pi_{\Ind_H^G U}^G X_R^{t_i}
\]
for the norm induced by $\mu_{X_R^{t_i}}$.
Write $\partial_i:\Omega X_R^{t_i}\to X_R$ and
$\partial_{i*}:\pi_U^K X_R^{t_i}\to\pi_{U-1}^K X_R$ for the $i$th boundary
of \eqref{eq:tate-compression} and its induced map, where
$K\in\{H,G\}$ and $U\in RO(K)$.

Let $R$ be $p$-torsionfree and $0\le i<h$.
Put $\lambda_j=\lambda_{p^h}(p^j)$ for $0\le j<h$.

We have the following definition.

\begin{definition}\label{def:intrinsic-orientations}
For $0\le j<h$, define
\[
\begin{aligned}
\mathcal O_{j,h}:(J_0^j)^{-1}B&\longrightarrow
\pi_{2-\lambda_j}^H X_R=\ker\Delta_{2-\lambda_j},\\
\ell&\longmapsto([\ell]_r)_{0\le r\le h}.
\end{aligned}
\]
Here $[\ell]_r$ is the class of $\ell$ in $A_r(2-\lambda_j)$.
Using $\Sigma_h$, the modules of representatives are
$(J_0^j)^{-1}B$ for $r\le j$ and $(J_0^{r-1})^{-1}B$ for $r>j$.
The inclusions $(J_0^j)^{-1}B\subseteq (J_0^{r-1})^{-1}B$ give
$\ell_r([\ell]_r)=\rho_r([\ell]_{r+1})$.
Under $\pi_\star^H X_R\to\pi_\star^H X_R^{t_i}$, a class
$\mathcal O_{j,h}(b')$ with $b'$ an $A$-basis element of $(J_0^j)^{-1}B$
maps to a unit for $i\le j<h$. The Euler class $a_{\lambda_j}$ maps to
a unit for $0\le j\le i$.

Choose $I_0=(\xi)$ and $b$ as in
Proposition~\ref{lem:period-repair-frame}. Put
\[
u_{\lambda_j}=\mathcal O_{j,h}\left(
\frac{b}{\prod_{l=1}^j\varphi^l(\xi)}\right).
\]
For $U=\rho_H(\mathbf a)$, $\mathbf a\in D_h^{(0)}$, the $A/J_i^h$-module
$\pi_U^H X_R^{t_i}=T_i(U)$ is free of rank one with basis element
\begin{equation}\label{eq:tate-unit-generator}
g_{T_i}^H(\mathbf a)=
\left(\frac{u_{\lambda_i}}{a_{\lambda_i}}\right)^{a_{i+1}}
\prod_{i<j<h}u_{\lambda_j}^{a_{j+1}-a_j}\prod_{0\le j\le i}a_{\lambda_j}^{a_{j+1}-a_j}.
\end{equation}
\end{definition}

Under $\Sigma_h^{a_h}$, this basis is represented by
\[
b^{a_h}\prod_{1\le r\le i}\varphi^r(\xi)^{-a_h}
\prod_{i<r\le h}\varphi^r(\xi)^{a_r-a_h}.
\]

For $s\in A/J_i^h$, choose a lift
$\widetilde s\in\pi_0^H X_R=A/(I_0J_0^h)$. Naturality gives
\[
N_H^G(s)=[N_H^G(\widetilde s)]\in A/J_i^n,
\]
where Proposition~\ref{prop:vertex-norms} with $U=0$ computes the
right-hand side. The result is independent of the lift.

It remains to compute the norm of the basis element $g_{T_i}^H(\mathbf a)$.

\begin{proposition}\label{lem:period-repair-norm}
Let $R$ be $p$-torsionfree, $0\le i<h$, and
$U=\rho_H(\mathbf a)$, $\mathbf a\in D_h^{(0)}$.
For $s\in A/J_i^h$, with $g_{T_i}^H(\mathbf a)$ as in
Definition~\ref{def:intrinsic-orientations}, we have
\[
\begin{aligned}
N_H^G(sg_{T_i}^H(\mathbf a))
={}&N_H^G(s)
\left(\frac{N_H^G(u_{\lambda_i})}
{a_{\Ind_H^G\lambda_i}}\right)^{a_{i+1}}\\
&\cdot\prod_{i<j<h}
N_H^G(u_{\lambda_j})^{a_{j+1}-a_j}
\prod_{0\le j\le i}
a_{\Ind_H^G\lambda_j}^{a_{j+1}-a_j}
\in T_i(\Ind_H^G U).
\end{aligned}
\]
\end{proposition}

\begin{proof}
Apply multiplicativity and
$N_H^G(a_{\lambda_j})=a_{\Ind_H^G\lambda_j}$
to \eqref{eq:tate-unit-generator}.
Use Proposition~\ref{prop:vertex-norms} for $N_H^G(u_{\lambda_j})$. Thus the proof follows.
\end{proof}

\subsubsection{Norms in shifted degrees}
Let $R$ be $p$-torsionfree and $U=\rho_H(\mathbf a)$, $\mathbf a\in D_h$.
Put $V=U-1$ and $T=\Ind_H^G V+1$.
By Proposition~\ref{prop:induction-of-degrees}, the degree coordinates of $T$ are
\begin{equation*}
\tfrac12\dim T^{C_{p^j}}=
\begin{cases}
p^{n-h}a_j-(p^{n-h}-1)/2,&0\le j\le h,\\
p^{n-j}a_h-(p^{n-j}-1)/2,&h<j\le n.
\end{cases}
\end{equation*}
For $n=h+1$, put
\begin{equation*}
J(\mathbf a)=(pa_0-(p-1)/2,\ldots,pa_h-(p-1)/2,a_h).
\end{equation*}

In degree $U-1$, we write a class by its edge coordinates, and then
compute the norm of each boundary class $\partial_{i*}z_i$.

\begin{proposition}\label{prop:edge-norms}
Let $R$ be $p$-torsionfree, $0\le i<h$, $U\in RO(H)$ be natural,
and $z\in\pi_U^H X_R^{t_i}$. Then
\begin{equation}\label{eq:full-boundary-norm}
N_H^G(\partial_{i*}z)
=\partial_{i*}\bigl(a_H^GN_H^G(z)\bigr).
\end{equation}
The norm is additive on $\pi_{U-1}^H X_R$.
For $U=\rho_H(\mathbf a)$, $\mathbf a\in D_h$, and $z_i\in T_i(U)$,
\[
N_H^G([z_0,\ldots,z_{h-1}])
=\sum_{i=0}^{h-1}
\partial_{i*}\bigl(a_H^GN_H^G(z_i)\bigr).
\]
\end{proposition}

\begin{proof}
The map of nonunital algebras $\partial_i$ in
Proposition~\ref{prop:boundary-indexed-multiplications}
gives \eqref{eq:full-boundary-norm}.

For $[G:H]=p$ and $x,y\in\pi_{U-1}^H X_R$, every term other
than $N_H^G(x)$ and $N_H^G(y)$ in
\cite[equation~(8.34)]{ZhangHZ} contains a product of
$p$ shifted classes, which is zero by \eqref{eq:path-cup}.
Hence
\[
N_H^G(x+y)=N_H^G(x)+N_H^G(y).
\]
Induction preserves shifted degrees, so norm transitivity
gives additivity for all $H<G$.
The formula for $N_H^G([z_0,\ldots,z_{h-1}])$ follows from
\eqref{eq:full-boundary-norm}
and
\[
[z_0,\ldots,z_{h-1}]
=\sum_{i=0}^{h-1}\partial_{i*}z_i.\qedhere
\]
\end{proof}

Let $R$ be $p$-torsionfree and $U=\rho_H(\mathbf a)$, $\mathbf a\in D_h$.
Use Proposition~\ref{lem:period-repair-norm} for $\mathbf a\in D_h^{(0)}$.
For $\mathbf a\in D_h^{(1)}$, the calculation is as follows.
For $U\notin RO_{\mathrm{sel}}(H)$, see Section~\ref{subsubsec:arbitrary-degree-norms}.

At $p=2$, let $h\ge1$ and put $m=2^{n-h}$. We have
\[
\Ind_H^G\sigma_{2^h}
\cong\bigoplus_{a=0}^{m/2-1}
\lambda_{2^n}((2a+1)2^{h-1}).
\]

Let $0\le i<h$ and $z\in\pi_U^H X_R^{t_i}$.
If $\epsilon(\mathbf a)=1$, the cofiber sequence
$
(H/C_{2^{h-1}})_+\longrightarrow S^0
\xrightarrow{a_{\sigma_{2^h}}}S^{\sigma_{2^h}}
$
gives
\[
a_{\sigma_{2^h}}:\pi_{U+\sigma_{2^h}}^H X_R^{t_i}
\twoheadrightarrow\pi_U^H X_R^{t_i}.
\]
Choose $y$ with $a_{\sigma_{2^h}}y=z$. Since
$U+\sigma_{2^h}=\rho_H(a_0+\tfrac12,\ldots,a_{h-1}+\tfrac12,a_h)
\in\rho_H(D_h^{(0)})$,
Proposition~\ref{lem:period-repair-norm} applies to $N_H^G(y)$, and
\[
N_H^G(z)
=N_H^G(a_{\sigma_{2^h}}y)
=a_{\Ind_H^G\sigma_{2^h}}N_H^G(y).
\]

This determines $N_H^G(z)$ for $\mathbf a\in D_h^{(1)}$. For $H=e$ and shifted $V$, $\pi_V^eX_R=0$.

\subsubsection{Norms in arbitrary representation degrees}
\label{subsubsec:arbitrary-degree-norms}
We now extend the norms from $RO_{\mathrm{sel}}(H)$ to $RO(H)$.
Let $H=C_{p^h}<G=C_{p^n}$ and write $V=V_0+D$ as in
Definition~\ref{def:character-translation}.
For $x\in\pi_V^H X_R$, multiplicativity gives, as in the proof of
\cite[Proposition~8.29]{ZhangHZ},
\[
N_H^G(x)=N_H^G(v_D)N_H^G(v_D^{-1}x),
\qquad v_D^{-1}x\in\pi_{V_0}^H X_R.
\]

Complex induction gives
\begin{equation*}
\Ind_H^G\lambda_{p^h}(k)
\cong\bigoplus_{a=0}^{p^{n-h}-1}
\lambda_{p^n}(k+ap^h).
\end{equation*}

For $b(t)=\min(t,p^n-t)$, the standard real weights give
\[
\widetilde v_{\Ind_H^G D}
=\prod_{r,j}\prod_{a=0}^{p^{n-h}-1}
\left(
\frac{\widetilde v_{b(rp^j+ap^h)/p^j,j}}
     {\widetilde v_{b(p^j+ap^h)/p^j,j}}
\right)^{m_{r,j}}.
\]

\begin{proposition}\label{lem:ordered-sign-units}
Let $R$ be $p$-torsionfree. Put
\[
q_l=
\frac{
\deg\bigl(\widetilde v_D^{C_{p^{\min(l,h)}}}\bigr)^{
p^{n-\max(l,h)}}}
{\deg\bigl(\widetilde v_{\Ind_H^G D}^{C_{p^l}}\bigr)}
\qquad(0\le l\le n).
\]
Then
\begin{equation*}
N_H^G(v_D)v_{\Ind_H^G D}^{-1}
=
(\widetilde\theta_{n+1}\varphi)^{-1}
\left(
q_n+\sum_{j=0}^{n-1}
\frac{q_j-q_{j+1}}{p^{n-j}}V^{n-j}(1)
\right).
\end{equation*}
\end{proposition}

\begin{proof}
Let $A(G)$ be the Burnside ring and put
\[
u=N_H^G(\widetilde v_D)\,
\widetilde v_{\Ind_H^G D}^{-1}\in A(G)_{(p)}^\times.
\]
Attention that for $k=rp^j,p^j$, folding the weight $k+ap^h$ requires conjugation
when $a\ge\lceil p^{n-h}/2\rceil$. For $l\le j$ the two signs
cancel. For $l>j$ both fixed spaces vanish. Thus the $C_{p^l}$-mark
of $u$ is $q_l$, and
\[
u=q_n[G/G]+\sum_{j=0}^{n-1}
\frac{q_j-q_{j+1}}{p^{n-j}}[G/C_{p^j}].
\]
Apply the unit map $A(G)_{(p)}\to\pi_0^G X_R$, which sends
\[
[G/C_{p^j}]\longmapsto\tr_{C_{p^j}}^G(1)
=(\widetilde\theta_{n+1}\varphi)^{-1}\bigl(V^{n-j}(1)\bigr).\qedhere
\]
\end{proof}

\subsection{Cyclotomic restriction and norms}\label{subsec:cyclotomic-maps}

We compute $R_{\mathrm{cyc}}$ on the kernel and cokernel in
\eqref{eq:uniform-groups}.
For perfect prisms, the degree-zero norm and cyclotomic restriction
are described in \cite[Section~3.1]{SulymaPrisms}. We start with the following definition.

\begin{definition}\label{def:cyclotomic-coefficients}
Let $R$ be as in Definition~\ref{def:natural-degree}, let $n\ge1$,
and let $U\in RO_{\mathrm{sel}}(C_{p^n})$ be natural. Put
\[
W=\varrho_p^*(U^{C_p}).
\]
For $p$-torsionfree $R$, define
\[
\begin{aligned}
\varphi^{-1}:A_{i+1}(U)\cong\varphi^*A_i(W)
&\longrightarrow A_i(W)
&& (0\le i<n),\\
\varphi^{-1}:T_{i+1}(U)\cong\varphi^*T_i(W)
&\longrightarrow T_i(W)
&& (0\le i<n-1),\\
a\otimes q&\longmapsto\varphi^{-1}(a)q,
\end{aligned}
\]
where the identifications use
\[
\varphi^*(\varphi^{(n-1)*}B)\cong\varphi^{n*}B,
\qquad
\varphi^*I_j\cong I_{j+1}.
\]
\end{definition}

We now describe cyclotomic restriction on vertex and edge coordinates
and compute its iterates on classes of the form $N_H^G(x)$.

\begin{proposition}\label{prop:cyclotomic-operations}
With $R,n,U,W$ as above, cyclotomic restriction is
\[
\begin{aligned}
R_{\mathrm{cyc}}:\ker\Delta_U
&\longrightarrow\ker\Delta_W,\\
(x_0,\ldots,x_n)
&\longmapsto
(\varphi^{-1}x_1,\ldots,\varphi^{-1}x_n),\\[4pt]
R_{\mathrm{cyc}}:\coker\Delta_U
&\longrightarrow\coker\Delta_W,\\
[z_0,\ldots,z_{n-1}]
&\longmapsto
[\varphi^{-1}z_1,\ldots,\varphi^{-1}z_{n-1}],
\end{aligned}
\]
under \eqref{eq:uniform-groups}.
The second map is zero for $n=1$.

For $V=V_0+D$ as in
Definition~\ref{def:character-translation} and
$x_0\in\pi_{V_0}^{C_{p^n}}X_R$,
\[
R_{\mathrm{cyc}}(v_Dx_0)
=
\left(\prod_r r^{m_{r,0}}\right)
v_{\varrho_p^*(D^{C_p})}R_{\mathrm{cyc}}(x_0).
\]
For $p=2$,
\[
R_{\mathrm{cyc}}(a_{\sigma_{2^n}})
=
\begin{cases}
a_{\sigma_{2^{n-1}}},&n\ge2,\\
1,&n=1.
\end{cases}
\]

For $H=C_{p^h}\le C_{p^n}$, $V\in RO(H)$,
$x\in\pi_V^H X_R$, and $0\le a\le n$,
\begin{equation}\label{eq:all-cyclotomic-crossings}
R_{\mathrm{cyc}}^aN_H^{C_{p^n}}(x)
=
\begin{cases}
N_{C_{p^{h-a}}}^{C_{p^{n-a}}}
R_{\mathrm{cyc}}^a(x),&a\le h,\\
N_e^{C_{p^{n-a}}}\psi^{p^{a-h}}
R_{\mathrm{cyc}}^h(x),&a>h.
\end{cases}
\end{equation}
Use the orbit diagonals and $\Phi^{C_{p^a}}\eta$ on the target.
\end{proposition}

\begin{proof}
Put $Y=S^{-U}\wedge X_R$.
Cyclotomic restriction is induced by
\[
Y^{C_{p^n}}
\longrightarrow
(Y^{\Phi_1})^{C_{p^n}}
\simeq (\Phi^{C_p}Y)^{C_{p^n}/C_p}
\xrightarrow[\mathrm{cyclotomic}]{\simeq}
(S^{-W}\wedge X_R)^{C_{p^{n-1}}}.
\]
By \eqref{eq:tate-compression}, $Y\to Y^{\Phi_1}$
projects to the diagram with $Y^{h_0}$ and $Y^{t_0}$ removed.

Assume $R$ is $p$-torsionfree. Since
$W^{C_{p^i}}\cong U^{C_{p^{i+1}}}$,
Proposition~\ref{prop:app-finite-comparison} and
\eqref{eq:app-sign-local-kernels} give
\[
\begin{aligned}
\varphi^{-(i+1)*}A_{i+1}(U)
&\cong
\pi_0\bigl((S^{-W^{C_{p^i}}}\wedge X_R)
^{hC_{p^{n-i-1}}}\bigr)
\cong\varphi^{-i*}A_i(W),
&&0\le i<n,\\
\varphi^{-(i+1)*}T_{i+1}(U)
&\cong
\pi_0\bigl((S^{-W^{C_{p^i}}}\wedge X_R)
^{tC_{p^{n-i-1}}}\bigr)
\cong\varphi^{-i*}T_i(W),
&&0\le i<n-1.
\end{aligned}
\]
The map $Y\to Y^{\Phi_1}$ induces the identity on each middle group.
By Subsection~\ref{subsec:period-line} and
\eqref{eq:app-moving-arrow}--\eqref{eq:app-frob-denominator},
for $0\le i<n-1$,
\[
\begin{aligned}
\ell_i\circ\varphi^{-1}
 &=\varphi^{-1}\circ\ell_{i+1}
 &&:A_{i+1}(U)\longrightarrow T_i(W),\\
\rho_i\circ\varphi^{-1}
 &=\varphi^{-1}\circ\rho_{i+1}
 &&:A_{i+2}(U)\longrightarrow T_i(W).
\end{aligned}
\]
For $\epsilon(\mathbf a)=1$, restrict to the kernels in
\eqref{eq:app-sign-local-kernels} by naturality of the
sign Euler cofiber sequence.
Thus we have a commutative diagram of abelian groups
\[
\begin{tikzcd}[column sep=large,row sep=large]
C_{C_{p^n}}^0(U)
  \arrow[r,"\Delta_U"]
  \arrow[d,"{(x_i)\mapsto(\varphi^{-1}x_{i+1})}"']
& C_{C_{p^n}}^1(U)
  \arrow[d,"{(z_i)\mapsto(\varphi^{-1}z_{i+1})}"] \\
C_{C_{p^{n-1}}}^0(W)
  \arrow[r,"\Delta_W"']
& C_{C_{p^{n-1}}}^1(W).
\end{tikzcd}
\]
Naturality of the homotopy fiber sequence gives the stated
maps on the kernel and cokernel under \eqref{eq:uniform-groups}.

For $0\le a\le h$, \cite[Theorem~2.35]{NormTHH} gives
\[
\Phi^{C_{p^a}}N_{C_{p^h}}^{C_{p^n}}
\simeq
N_{C_{p^{h-a}}}^{C_{p^{n-a}}}\Phi^{C_{p^a}}.
\]
Together with the norm compatibility of the cyclotomic equivalence,
this gives the first case of \eqref{eq:all-cyclotomic-crossings}.
By \cite[Remark~5.2]{SulymaNew} and norm transitivity,
\[
R_{\mathrm{cyc}}N_e^{C_{p^m}}
=
N_e^{C_{p^{m-1}}}R_{\mathrm{cyc}}N_e^{C_p}
=
N_e^{C_{p^{m-1}}}\psi^p
\qquad(m\ge1).
\]
\end{proof}

\appendix
\section{Homotopy fixed points and Tate spectra}
\label{app:spectral-foundations}

In this appendix we compute the homotopy groups
of \eqref{eq:tate-compression}.
Let $R$ be a nonzero $p$-torsionfree perfectoid ring, with
$A,I_j,J_0^r,B,\Phi_B$ as in Section~\ref{sec:uniform}.
Every $I_j$ is a Cartier ideal.

\subsection{Homotopy groups of homotopy fixed-point and Tate spectra}

We first treat degrees that lift to virtual complex representations of
the circle. We use the same lift at each finite subgroup to pass from
the circle calculation to homotopy fixed points and Tate spectra for $C_{p^r}$.

\begin{definition}\label{def:app-local-spectra}
Let $U=\rho_p(\mathbf a)$ with $\mathbf a\in D_n^{(0)}$. The map
$\Sigma_n^{a_n}$ gives the identification
\[
E_0(U)\xrightarrow{\sim}B^{a_n}
        \prod_{j=1}^{n-1}(J_0^j)^{a_j-a_{j+1}}.
\]

Let $\mathbb C(k)$ be the complex $\mathbb T$-representation
of weight $k$. Put
\[
\widetilde U
=a_n\mathbb C(0)
+\sum_{j=0}^{n-1}(a_j-a_{j+1})\mathbb C(p^j).
\]
Its underlying real representation restricts to $U$.
For every $r$, we use the restriction of this same
$\widetilde U$ to $C_{p^r}$ to compute the homotopy groups.
Thus $E_0(U)$ is defined using $\widetilde U$, with $n$ fixed.

\end{definition}

By \cite[Section~2, p.~12]{SulymaNew},
$\THH(R\mid\mathbb S_{W(R^\flat)};\Z_p)\simeq X_R$
in $\Sp^{B\mathbb T}$.
Apply \cite[Theorem~1.9 and Remarks~3.1,~3.5]{SulymaNew}
with sequence $(a_0,\ldots,a_n;a_n)$.
After pullback along $\varphi$, using
$\varphi^*\operatorname{BK}(A,I_0)\cong B$, we obtain
\[
\begin{aligned}
\pi_0\bigl((S^{-\widetilde U}\wedge X_R)^{h\mathbb T}\bigr)
&\cong I_0^{a_0^+}E_0(U),\\
\pi_0\bigl((S^{-\widetilde U}\wedge X_R)^{t\mathbb T}\bigr)
&\cong E_0(U).
\end{aligned}
\]
These spectra have no odd homotopy groups, and the canonical map
is induced by $I_0^{a_0^+}E_0(U)\hookrightarrow E_0(U)$.

On restriction to $C_{p^r}$, the canonical map is still induced by this
inclusion, now on the quotients described below.

\begin{proposition}\label{prop:app-finite-comparison}
Let $U=\rho_p(\mathbf a)$, $\mathbf a\in D_n^{(0)}$,
with the $\mathbb T$-lift of Definition~\ref{def:app-local-spectra}.
For $r\geq0$ and $k\in\mathbb Z$,
\begin{align}
 \pi_0\bigl((S^{-U}\wedge X_R)^{hC_{p^r}}\bigr)&\cong
 \frac{I_0^{a_0^+}E_0(U)}{J_0^r I_0^{(a_0+1)^+}E_0(U)},
 &\pi_{2k+1}\bigl((S^{-U}\wedge X_R)^{hC_{p^r}}\bigr)&=0,
 \label{eq:app-finite-borel}\\
 \pi_0\bigl((S^{-U}\wedge X_R)^{tC_{p^r}}\bigr)&\cong \frac{E_0(U)}{J_0^rE_0(U)},
 &\pi_{2k+1}\bigl((S^{-U}\wedge X_R)^{tC_{p^r}}\bigr)&=0.
 \label{eq:app-finite-tate}
\end{align}
The canonical map from homotopy fixed points to the Tate spectrum is induced by
$I_0^{a_0^+}E_0(U)\hookrightarrow E_0(U)$.
\end{proposition}

\begin{proof}
Locally choose $I_0=(\xi)$ and $B=Ae$, with dual basis
$e^\vee\in B^\vee=\pi_{-2}X_R^{h\mathbb T}$.
By \cite[Lemma~3.1 and its proof]{Riggenbach}, for each $m\ge0$
we have the cofiber sequence
\[
\Sigma^{-2}M_m^{h\mathbb T}
\xrightarrow{\left(\prod_{j=1}^{r}\varphi^j(\xi)\right)e^\vee}
M_m^{h\mathbb T}
\xrightarrow{\operatorname{res}}
M_m^{hC_{p^r}},
\qquad
M_m=S^{m\mathbb C(1)-\widetilde U}\wedge X_R.
\]

For $K=\mathbb T,C_{p^r}$, the model
$EK=S(\mathbb C(1)^{\oplus\infty})|_K$ gives
\[
\operatorname*{colim}_{m\ge0}M_m^{hK}
\simeq
\bigl(S^{-\widetilde U}\wedge X_R\bigr)^{tK},
\]
where the transition maps are induced by
$S^0\to S^{\mathbb C(1)}$.
Under this equivalence, the map from $m=0$ is
$\operatorname{can}$.
Taking colimits of the cofiber sequences above gives the
bottom row of the following commutative diagram;
the top row is the sequence at $m=0$:
\[
\begin{tikzcd}[column sep=large,row sep=huge]
\Sigma^{-2}(S^{-\widetilde U}\wedge X_R)^{h\mathbb T}
  \arrow[r,"{\left(\prod_{j=1}^r\varphi^j(\xi)\right)e^\vee}"]
  \arrow[d,"\Sigma^{-2}\operatorname{can}"']
&
(S^{-\widetilde U}\wedge X_R)^{h\mathbb T}
  \arrow[r,"\res"]
  \arrow[d,"\operatorname{can}"]
&
(S^{-\widetilde U}\wedge X_R)^{hC_{p^r}}
  \arrow[d,"\operatorname{can}"]
\\
\Sigma^{-2}(S^{-\widetilde U}\wedge X_R)^{t\mathbb T}
  \arrow[r,"{\left(\prod_{j=1}^r\varphi^j(\xi)\right)
    \operatorname{can}(e^\vee)}"']
&
(S^{-\widetilde U}\wedge X_R)^{t\mathbb T}
  \arrow[r]
&
(S^{-\widetilde U}\wedge X_R)^{tC_{p^r}}.
\end{tikzcd}
\]

For $k\in\mathbb Z$, the left square on $\pi_{2k}$ is
\[
\begin{tikzcd}[column sep=large,row sep=large]
I_0^{(a_0+k+1)^+}B^{k+1}\otimes_A E_0(U)
  \arrow[r,hook]
  \arrow[d,hook]
&
I_0^{(a_0+k)^+}B^k\otimes_A E_0(U)
  \arrow[d,hook]
\\
B^{k+1}\otimes_A E_0(U)
  \arrow[r,hook]
&
B^k\otimes_A E_0(U),
\end{tikzcd}
\]
with horizontal maps
\[
ae^{k+1}\otimes\ell
\longmapsto
\left(\prod_{j=1}^r\varphi^j(\xi)\right)
ae^k\otimes\ell.
\]
Since
\[
J_0^r=\left(\prod_{j=1}^r\varphi^j(\xi)\right),
\]
the long exact sequences give
\eqref{eq:app-finite-borel}--\eqref{eq:app-finite-tate}.
The right vertical inclusions induce
\[
\operatorname{can}:
\frac{I_0^{a_0^+}E_0(U)}
     {J_0^rI_0^{(a_0+1)^+}E_0(U)}
\longrightarrow
\frac{E_0(U)}{J_0^rE_0(U)},
\qquad
[x]\longmapsto[x].
\]
\end{proof}

\subsection{Maps from geometric fixed points to Tate spectra}

We next compute the map from geometric fixed points to Tate spectra.
For $r\ge1$ and $U\in RO(C_{p^r})$, let
\[
\theta_U:
\Phi^{C_p}(S^{-U}\wedge X_R)
\longrightarrow (S^{-U}\wedge X_R)^{tC_p}
\]
be the canonical map. Let
$q:C_{p^r}\to C_{p^r}/C_p$ and write
$U=q^*W+M$, with $W\in RO(C_{p^r}/C_p)$ and $M^{C_p}=0$.
For actual $M$, the sphere $S(M)$ is $C_p$-free, so the Euler map
$a_M$ is a Tate equivalence and $\Phi^{C_p}(a_M)=\mathrm{id}$.
For virtual $M$, it's implied by inverse Euler equivalences for negative summands.
By naturality in \cite[Lemma~II.4.5]{NS}, we have
\begin{equation}\label{eq:app-moving-arrow}
\theta_U=(a_M^{tC_p})^{-1}\circ(S^{-W}\wedge\theta_0).
\end{equation}
Take $r=n\ge1$ and
$U=\rho_p(\mathbf a)$ with $\mathbf a\in D_n^{(0)}$.
Under the cyclotomic equivalence, $\theta_0$ is cyclotomic
Frobenius. Put
\[
W=\rho_{C_{p^{n-1}}}(a_1,\ldots,a_n),
\qquad U=q^*W+(a_0-a_1)\lambda_0.
\]
Definitions~\ref{def:natural-degree}
and~\ref{def:app-local-spectra} give
\begin{align}
\varphi^*E_0(W)
  &=I_1^{-a_1}E_0(U),\notag\\
\varphi^*\bigl(I_0^{a_1^+}E_0(W)\bigr)
  &=I_1^{(-a_1)^+}E_0(U).
\label{eq:app-frob-numerator}
\end{align}

For the circle lift of Definition~\ref{def:app-local-spectra},
the Frobenius arrow in \cite[Theorem~1.9, Notation~3.14 and Example~3.15]{SulymaNew},
using the product of fractional ideals for virtual weights,
with the Breuil--Kisin Frobenius of Subsection~\ref{subsec:period-line},
has $A$-linearization
\[
I_1^{(-a_1)^+}E_0(U)\hookrightarrow E_0(U)
\]
under \eqref{eq:app-frob-numerator}.

By \cite[Lemma~II.4.1]{NS}, restriction to
$C_{p^{n-1}}\le\mathbb T/C_p$ gives
\[
(S^{-W}\wedge X_R)^{hC_{p^{n-1}}}
\xrightarrow{\ \theta_U^{hC_{p^{n-1}}}\ }
(S^{-U}\wedge X_R)^{tC_{p^n}}.
\]
Since $\varphi(J_0^{n-1})=I_1^{-1}J_0^n$,
\begin{equation}\label{eq:app-frob-denominator}
\varphi^*\bigl(J_0^{n-1}I_0^{(a_1+1)^+}E_0(W)\bigr)
=J_0^n I_1^{(-a_1-1)^+}E_0(U)
\subseteq J_0^nE_0(U).
\end{equation}
Naturality and Proposition~\ref{prop:app-finite-comparison}
identify the induced finite quotient map with $[x]\mapsto[x]$.
Applying this calculation to $U^{C_{p^i}}$ at $C_{p^{n-i}}$
and pulling back along $\varphi^i$ gives $\ell_i$ and $\rho_i$
as the maps $[x]\mapsto[x]$ of Definition~\ref{def:natural-degree}.
For $i=n-1$, Proposition~\ref{prop:app-finite-comparison}
with $r=0$ gives
$\varphi^{n*}\pi_{2a_n}X_R=A_n(U)$.

\subsection{The sign representation at \texorpdfstring{$p=2$}{p=2}}

What remains is the selected natural degrees with an odd
coefficient of the sign representation. By the sign Euler cofiber sequence, it suffices to compute the local terms' kernels of restriction in complex degrees.
Assume $p=2$ and $n\geq1$. Let $U$ be a
selected natural degree with $\epsilon(\mathbf a)=1$, and put
$W=U-\sigma_{2^n}$. Then $W$ is complex, with degree coordinates
$(a_0-\tfrac12,\ldots,a_{n-1}-\tfrac12,a_n)$.
For $0\le i<n$, put $r=n-i$. Dualizing
\[
(C_{2^r}/C_{2^{r-1}})_+\longrightarrow S^0
\longrightarrow S^{\sigma_{2^r}}
\]
and applying Shapiro identifies the local terms for $U^{C_{2^i}}$
with the fibers of restriction from $C_{2^r}$ to $C_{2^{r-1}}$
for $W^{C_{2^i}}$, for both homotopy fixed points and Tate spectra.
Use the same circle lift at both subgroups.
By Proposition~\ref{prop:app-finite-comparison}, these restrictions
are surjective quotient maps in even degrees, so the fibers have
no odd homotopy groups.
In degree zero, we have
\begin{align}
 \pi_0\bigl((S^{-U^{C_{2^i}}}\wedge X_R)^{hC_{2^r}}\bigr)
   &\cong\frac{J_0^{r-1}I_0^{(a_i+1/2)^+}E_0(W^{C_{2^i}})}
           {J_0^rI_0^{(a_i+1/2)^+}E_0(W^{C_{2^i}})},\notag\\
 \pi_0\bigl((S^{-U^{C_{2^i}}}\wedge X_R)^{tC_{2^r}}\bigr)
   &\cong\frac{J_0^{r-1}E_0(W^{C_{2^i}})}{J_0^rE_0(W^{C_{2^i}})}.
 \label{eq:app-sign-local-kernels}
\end{align}

After pulling back along $\varphi^i$,
\[
J_i^{n-1}E_i(W)=\mathcal S_i(U),
\]
so \eqref{eq:app-sign-local-kernels} gives \eqref{eq:uniform-sign}.

Restrict $\ell_i,\rho_i$ for $W$ to these kernels.
If $r\ge2$, then $C_2\subseteq\ker(\sigma_{2^r})$, so
$(\sigma_{2^r})^{C_2}$ is the sign representation of $C_{2^r}/C_2$.
If $r=1$, then $(\sigma_{2^r})^{C_2}=0$ and $\Phi^{C_2}(a_{\sigma_2})$
is the identity. The map $\rho_{n-1}$ follows by naturality
of the sign Euler cofiber sequence.

\section{Proofs for norms}
\label{app:operation-foundations}

For a nonzero $p$-torsionfree perfectoid ring $R$, we prove two formulas
used in Subsection~\ref{subsec:graded-norms}.

\subsection{Maps induced by geometric fixed points}
\label{app:ghosts}

We compute the maps $\gamma_i$ of
Definition~\ref{def:geometric-fixed-point-maps}
for $U=\rho_p(\mathbf a)$ with $\mathbf a\in D_n^{(0)}$.

Put $L=I_0/I_0^2$, $\lambda=\mathbb C(1)$. Identify $\pi_{2a}X_R=L^a$ for $a\ge0$
by the restriction $L\xrightarrow{\sim}\pi_2X_R$ in degree $\lambda$
\cite[Example~3.2 and Remark~3.5]{SulymaNew}.
For $P=X_R,\varrho_p^*(X_R^{tC_p})$, the Euler map on the two-cell of
$B\mathbb T$ sends a weight-one class restricting to $z$ to
$c_1(\lambda)z\in H^2(B\mathbb T;\pi_2P)$.
With the same generator via $\varrho_p$, cyclotomic Frobenius induces
$\varphi:L\to I_1/I_1^2$ on these quotients and
$x\mapsto\Phi_B(1\otimes x)$ on $\pi_2$ of circle homotopy fixed points
(Subsection~\ref{subsec:period-line}, \cite[Proposition~6.2 and Remark~6.6]{BMS}).
Hence integer-degree restriction is
\[
\delta:I_0B/I_0^2B\xrightarrow{\sim}L,\qquad
[x]\longmapsto\varphi^{-1}\!\left(\operatorname{pr}_1\Phi_B(1\otimes x)\right),
\]
where $\operatorname{pr}_1:B\to B/I_1B=I_1/I_1^2$ is the first projection.
Write
\begin{equation}\label{eq:geometric-kappa}
\delta=\mathrm{id}_L\otimes\kappa,\qquad
\kappa:B/I_0B\xrightarrow{\sim}R.
\end{equation}
The change to HKR coordinates is $\delta^{-1}$ followed by the second
map in \eqref{eq:period-cotangent-reduction}.

\begin{definition}\label{def:geometric-fixed-point-coefficients}
For $0\le i<n$ and $a_i\ge0$, use $\Sigma_{n-i}$ from \eqref{eq:period-repair-data} to identify the lines
\[
\varphi^{-i*}\bigl(I_i^{a_i}E_i(U)\bigr)
\cong I_0^{a_i}B^{a_n}
\prod_{i<j<n}(J_0^{j-i})^{a_j-a_{j+1}}.
\]
Define
\[
\bar\theta_{a_i}:
\bigl(\varphi^{-i*}\bigl(I_i^{a_i}E_i(U)/I_i^{a_i+1}E_i(U)\bigr)\bigr)[1/p]
\longrightarrow (I_0/I_0^2)^{a_i}[1/p]
\]
by tensoring the identity on $(I_0/I_0^2)^{a_i}$ with the
$a_n$th tensor power of $\kappa$ from \eqref{eq:geometric-kappa} and the
$(a_j-a_{j+1})$th tensor powers of
\[
J_0^{j-i}/I_0J_0^{j-i}\longrightarrow R,
\qquad [a]\longmapsto\theta(a)
\quad(i<j<n).
\]
For negative powers, use the duals of the inverse maps after
inverting $p$.
\end{definition}

\begin{proposition}\label{prop:geometric-fixed-point-coordinates}
Let $x=(x_0,\ldots,x_n)\in\ker\Delta_U$.
For $0\le i<n$, if $a_i<0$, then $\gamma_i(x)=0$.
If $a_i\ge0$, let $\overline{x_i}$ be the image of $x_i$ in
$I_i^{a_i}E_i(U)/I_i^{a_i+1}E_i(U)$.
Then
\begin{equation}\label{eq:ghost-line-conversion}
\gamma_i(x)
=p^{-\sum_{i<j<n}(a_j-a_n)}
\bar\theta_{a_i}\bigl(\varphi^{-i*}\overline{x_i}\bigr)
\in L^{a_i}.
\end{equation}
\[
p^{-\sum_{i<j<n}(a_j-a_n)}\bar\theta_{a_i}:
\varphi^{-i*}\bigl(I_i^{a_i}E_i(U)/I_i^{a_i+1}E_i(U)\bigr)
\xrightarrow[\ R\text{-linear}\ ]{\sim}(I_0/I_0^2)^{a_i}.
\]
At the terminal vertex,
\[
\gamma_n(x)=\delta^{\otimes a_n}\bigl(\varphi^{-n}(x_n)\bigr)
\in L^{a_n}\qquad(a_n\ge0),
\]
where $A_n(U)=\varphi^{n*}\bigl((I_0B/I_0^2B)^{\otimes_R a_n}\bigr)$.
If $a_n<0$, then $\gamma_n(x)=0$.
\end{proposition}

\begin{proof}
By \cite[Theorem~6.1]{BMS}, it suffices to consider $a_i\ge0$.
By \cite[Lemma~2.2.8]{BhattLurie} and $p$-torsionfreeness of $R$,
\[
J_0^k/I_0J_0^k\xrightarrow[\sim]{\ \theta/p^k\ }R
\quad(k\ge1).
\]
Taking tensor powers, with
$\sum_{i<j<n}(j-i)(a_j-a_{j+1})=\sum_{i<j<n}(a_j-a_n)$,
and then using $\kappa$ proves the isomorphism.

For $k\ge1$, pullback along $z\mapsto z^{p^k}$ induces
$I_0J_0^k\hookrightarrow I_0$ in the Euler coordinates
of \cite[Notation~2.11 and Example~3.2]{SulymaNew}. Its underlying
degree is $p^k$, so $p^k\res_e^{\mathbb T}(z)=[z]\in L$ for $z\in I_0J_0^k$.
Thus restriction in degrees $\mathbb C(p^k)-\lambda$ and
$\mathbb C(0)-\lambda$ is $\theta/p^k$ and $\kappa$, respectively.
Tensor products and duals give the restriction for $U^{C_{p^i}}$
\cite[Remark~3.1 and Examples~3.3--3.4]{SulymaNew}.

For the circle lift of Definition~\ref{def:app-local-spectra},
Proposition~\ref{prop:app-finite-comparison} identifies restriction
to $C_{p^{n-i}}$ with the quotient onto $\varphi^{-i*}A_i(U)$.
Since $\res_e^{\mathbb T}
=\res_e^{C_{p^{n-i}}}\res_{C_{p^{n-i}}}^{\mathbb T}$,
composition with the $i$th vertex of \eqref{eq:tate-compression}
gives \eqref{eq:ghost-line-conversion}.
\end{proof}

\subsection{Proof of \eqref{eq:full-boundary-norm}}
\label{app:full-cut-norms}

For $H=C_{p^h}<G=C_{p^n}$ and $0\le i<h$, we compute the norm of a
boundary class. We show that $\partial_i$ respects the indexed
multiplication for $G/H$.

\begin{proposition}\label{prop:boundary-indexed-multiplications}
Let $H=C_{p^h}<G=C_{p^n}$ and $0\le i<h$.
For the indexing system~\cite[Definition~3.22]{BHOperadic} generated by $G/H$, the Tate compression
diagram truncated at $X_R^{\Phi_{i+1}}$ is a diagram of algebras, and
\[
\partial_i:\Omega X_R^{t_i}\longrightarrow X_R
\]
is a map of nonunital algebras.
Consequently, \eqref{eq:full-boundary-norm} holds.
\end{proposition}

\begin{proof}
For $-1\le r\le i$, we have $C_{p^{r+1}}\le H$.
By \cite[Theorem~6.1 and Corollary~6.5]{BHNorm}, the functors
$\widetilde E\mathcal F_r\wedge(-)$ preserve these algebras,
and the localization maps are maps of algebras.
The $K$-equivariant diagonals $E\mathcal F_r\to\prod_T E\mathcal F_r$,
for all admissible finite $K$-sets $T$ and $K\le G$,
make $F((E\mathcal F_r)_+,-)$ a functor on these algebras,
and the completion maps are maps of algebras.

By Proposition~\ref{prop:compression}, the $\mathcal F_r$-completion
maps of $X_R^{t_{r-1}}$ are equivalences of these algebras
for $1\le r\le i$.
Naturality of completion gives the factorizations of
Proposition~\ref{prop:compression} in the same category.
Iterated isotropy separation gives
\begin{equation}\label{eq:full-cut-square}
X_R\simeq
\holim\bigl(X_R^{h_0}\to X_R^{t_0}\leftarrow\cdots\to
X_R^{t_{i-1}}\leftarrow X_R^{h_i}\bigr)
\mathop{\times}_{X_R^{t_i}}X_R^{\Phi_{i+1}},
\end{equation}
where the homotopy limit is $X_R^{h_0}$ for $i=0$.

For a homotopy pullback $B\times_D C$ of these algebras,
the zero maps to $B$ and $C$ give a map of nonunital algebras
\[
\Omega D=0\mathop{\times}_D0\longrightarrow B\times_D C.
\]
Limits of nonunital algebras are computed on spectra,
so the map for \eqref{eq:full-cut-square} is $\partial_i$.
The indexed multiplication on $\Omega D\simeq S^{-1}\wedge D$ is
\[
S^{-\mathbb R[G/H]}\wedge N_H^G\res_H^G D
\longrightarrow S^{-\mathbb R[G/H]}\wedge D
\xrightarrow{\ a_H^G\ }S^{-1}\wedge D.
\]
$a_H^G$ is dual to the one-point compactification of
\[
\mathbb R\hookrightarrow\mathbb R[G/H],
\qquad t\longmapsto(t)_{gH\in G/H}.
\]
Then taking homotopy groups gives \eqref{eq:full-boundary-norm}
(cf.\ \cite[Theorem~3.1.1]{Balderrama} for the Euler factor).
\end{proof}

\section{The computation for perfect \texorpdfstring{$\mathbb F_p$}{Fp}-algebras}
\label{app:characteristic-p}
In characteristic $p$, we obtain the graded Tambara functor from that
of $\HZ$ by inverting Euler classes and extending to Witt
vectors. We also compute cyclotomic restriction.

Let $k$ be a perfect $\F_p$-algebra, $A=W(k)$, and $\varphi$ its Witt
Frobenius.
Write $\underline A$ for the constant Tambara functor.
By \cite[Proposition~2.1]{KrauseNikolaus} and \cite{Bokstedt},
$\pi_*\THH(Hk)=k[\beta]$ with $|\beta|=2$.
Thus the underlying spectrum is $p$-complete, and the cyclotomic equivalences
in Section~\ref{sec:introduction} give $\THH(Hk)\simeq X_k$ on every $C_{p^n}$.

\subsection{Geometric fixed points and the cyclotomic restriction tower}

We identify $X_k$ with the geometric $C_p$-fixed points of
$H\underline A$ by taking the limit along cyclotomic restriction.
For this comparison to preserve norms, we take the limit in genuine
commutative rings.

For a genuine commutative $C_{p^{r+1}}$-ring $Y$, the map
\begin{equation*}
\omega_Y:Y^{C_p}\longrightarrow\Phi^{C_p}Y
\end{equation*}
preserves norms by
\cite[Theorem~6.1 and Corollary~6.5]{BHNorm} and
\cite[Proposition~2.19 and Theorem~2.35]{NormTHH}.
By \cite[Proposition~1.2, Theorem~2.36 and
Theorems~4.6--4.7]{NormTHH}, the cyclotomic equivalence
$\theta:\varrho_p^*\Phi^{C_p}X_k\xrightarrow{\sim}X_k$
also preserves norms, and
\[
R_{\mathrm{cyc}}=\theta\circ\varrho_p^*(\omega_{X_k}):
\varrho_p^*(X_k^{C_p})\longrightarrow X_k.
\]
For finite cyclic groups, we omit the pullbacks along $\varrho_p$.

\begin{definition}
For $r\ge0$, take the following limit in genuine commutative
$C_{p^r}$-rings:
\[
\mathcal R_r=\holim_{j\ge0,\,R_{\mathrm{cyc}}^{C_{p^j}}}
(X_k^{C_{p^j}})|_{C_{p^r}}.
\]
Let $p_r:\mathcal R_r\to X_k|_{C_{p^r}}$ be the projection.
The underlying limit is computed in genuine spectra.
\end{definition}

We start with the following observation.

\begin{proposition}\label{prop:cp-restriction-limit}
The map
\begin{equation}\label{cp-limit-geometric}
\theta\circ\Phi^{C_p}p_{r+1}:
\Phi^{C_p}\mathcal R_{r+1}\xrightarrow{\sim}X_k|_{C_{p^r}}
\end{equation}
is an equivalence of genuine commutative rings. There is an equivalence
$\iota_r:H\underline A\xrightarrow{\sim}\mathcal R_r$ inducing
\begin{equation}\label{cp-witt-gauge}
A\longrightarrow\pi_0^{C_{p^i}}\mathcal R_r=W(k),
\qquad a\longmapsto\varphi^{1-i}(a)
\quad(0\le i\le r).
\end{equation}
\end{proposition}

\begin{proof}
Since $X_k$ is connective, \eqref{cp-limit-geometric} is an
equivalence by \cite[Proposition~5.2, the proof of Corollary~5.4,
and Corollary~5.13]{AN}.
By \cite[Theorem~6.14]{AN},
\[
\pi_s^K\mathcal R_r=
\begin{cases}W(k),&s=0,\\0,&s\ne0,\end{cases}
\qquad K\le C_{p^r},
\]
and $\res_{C_{p^i}}^{C_{p^l}}=\varphi^{l-i}$.
The residual circle acts trivially on $W(k)$.
Since restriction is invertible, restriction determines transfers
and norms. Thus \eqref{cp-witt-gauge} is an isomorphism of
Tambara functors, which induces $\iota_r$ by
\cite[Theorems~5.2--5.3]{Ullman}.
\end{proof}

The comparison in the derived category of Mackey modules is given in
\cite[Theorem~3]{Kriz}. Using $\iota_{n+1}$, we obtain a comparison of
genuine commutative $C_{p^n}$-rings.

\begin{proposition}
\label{prop:cp-normed-comparison}
For every perfect $k$ and $n\geq0$, the composite
\[
c_{n,k}:\Phi^{C_p}H\underline A_{C_{p^{n+1}}}
 \xrightarrow{\Phi^{C_p}\iota_{n+1}}\Phi^{C_p}\mathcal R_{n+1}
 \xrightarrow{\Phi^{C_p}p_{n+1}}\Phi^{C_p}X_k
 \xrightarrow{\theta}X_k|_{C_{p^n}}
\]
is an equivalence of genuine commutative $C_{p^n}$-rings, compatible
with restriction to subgroups in the homotopy category. The map on $\pi_0^{C_{p^i}}$ is
\begin{equation}\label{cp-finite-gauge}
A/p^{i+1}\longrightarrow W_{i+1}(k),\qquad
 [a]\longmapsto[\varphi^{-i}(a)].
\end{equation}
For $k=\F_p$ the source can equally be
$\Phi^{C_p}H\underline\Z$.
\end{proposition}

\begin{proof}
The equivalence follows from \eqref{cp-limit-geometric}. At the preimage
$C_{p^{i+1}}$ of $C_{p^i}$, the map on degree-zero homotopy is
$\varphi^{-i}$. Projection followed by cyclotomic
restriction is Witt truncation, so we obtain
\eqref{cp-finite-gauge}.

For $k=\F_p$, the map
$\Phi^{C_p}H\underline\Z\to\Phi^{C_p}H\underline{\Z_p}$
is a genuine equivalence by \cite[Section~1, p.~3]{ZhangHZ}.
\end{proof}

We can therefore pass from $\F_p$ to $k$ on representation-graded
homotopy groups by extending scalars from $\Z_p$ to $W(k)$.

\begin{proposition}\label{prop:cp-flat-extension}
For $V\in RO(C_{p^n})$,
\begin{equation}\label{cp-flat-extension}
 \pi_V^{C_{p^n}}X_k\cong
A\otimes_{\Z_p}\pi_V^{C_{p^n}}X_{\F_p}.
\end{equation}
This identification preserves products, restrictions and transfers.
\end{proposition}
\begin{proof}
Since $A$ is flat over $\Z_p$, write
\[
A\cong\operatorname*{colim}_{\alpha}\Z_p^{r_\alpha}
\]
as a filtered colimit of finite free modules. Thus, as genuine
$C_{p^{n+1}}$-spectra,
\[
H\underline A\simeq
\operatorname*{colim}_{\alpha}
\bigl(H\underline{\Z_p}\bigr)^{\oplus r_\alpha}.
\]
Geometric fixed points preserve this colimit. Attention that representation
spheres are compact. Proposition~\ref{prop:cp-normed-comparison}
therefore gives
\[
\begin{aligned}
\pi_V^{C_{p^n}}X_k
&\cong
\pi_V^{C_{p^n}}\Phi^{C_p}H\underline A\\
&\cong
\operatorname*{colim}_{\alpha}
\left(
\pi_V^{C_{p^n}}\Phi^{C_p}H\underline{\Z_p}
\right)^{\oplus r_\alpha}\\
&\cong
\operatorname*{colim}_{\alpha}
\left(
\Z_p^{r_\alpha}\otimes_{\Z_p}
\pi_V^{C_{p^n}}X_{\F_p}
\right)\\
&\cong
A\otimes_{\Z_p}\pi_V^{C_{p^n}}X_{\F_p}.
\end{aligned}
\]
The inverse sends $a\otimes x$ to $ax$ and is induced by
$H\underline{\Z_p}\to H\underline A$ and Proposition~\ref{prop:cp-normed-comparison},
with the $A$-action given by \eqref{cp-finite-gauge}.
It preserves products and bilinearity, restrictions by naturality, and transfers by
Frobenius reciprocity.
\end{proof}

\subsection{Derivation from $H\underline{\mathbb{Z}}$}

We use Euler localization to derive the Tambara functor for $X_k$ from $\HZ$ and then compute cyclotomic restriction. We begin by determining how many powers of the Euler
class suffice in a given representation degree.

\begin{proposition}\label{lem:cp-stabilization}
Let $q:C_{p^{n+1}}\twoheadrightarrow C_{p^n}$ be the quotient map, and
let $\lambda=\lambda_{p^{n+1}}(1)$ be the faithful oriented plane.
For $H\leq C_{p^n}$ write $\widetilde H=q^{-1}H$ and
$\lambda_H=\lambda|_{\widetilde H}$.  For every $V\in\RO(H)$,
\begin{equation}\label{cp-euler-bound}
 \pi_{q^*V-m\lambda_H}^{\widetilde H}H\underline A
       \xrightarrow{\ \sim\ }\pi_V^H X_k,
 \qquad m\geq0,\quad 2m>\dim V.
\end{equation}
The transition $m\mapsto m+1$ is multiplication by $a_{\lambda_H}$.
\end{proposition}
\begin{proof}
Since $A$ is flat over $\Z$, the argument for
\eqref{cp-flat-extension} applies over $\Z$.
For $m\ge0$ and $2m>\dim V$,
\cite[Section~1, p.~3]{ZhangHZ} and
Proposition~\ref{prop:cp-normed-comparison} give
\[
\begin{aligned}
\pi_{q^*V-m\lambda_H}^{\widetilde H}H\underline A
&\cong
A\otimes_{\Z}
\pi_{q^*V-m\lambda_H}^{\widetilde H}\HZ\\
&\xrightarrow{\sim}
A\otimes_{\Z}\pi_V^H\Phi^{C_p}\HZ\\
&\cong
\pi_V^H\Phi^{C_p}H\underline A\\
&\xrightarrow{\sim}
\pi_V^H X_k.
\end{aligned}
\]
Under the first isomorphism,
$1\otimes(a_{\lambda_H}\cdot-)$ corresponds to
$a_{\lambda_H}\cdot-$.
\end{proof}

We next describe how to take the norm of a class represented at one of
these finite stages.

\begin{proposition}\label{prop:cp-induced}
The $RO(C_{p^n})$-graded Tambara functor of $X_k$ is obtained
from the $RO(C_{p^{n+1}})$-graded Tambara functor of $\HZ$ by localization of Euler classes and extension of
coefficients to $A$.
\end{proposition}

\begin{proof}
By Propositions~\ref{prop:cp-normed-comparison}
and~\ref{prop:cp-flat-extension} and naturality, we get the Green functor structure.
It remains to determine the norms.

Use the notation of Proposition~\ref{lem:cp-stabilization}.
Let $H\le L\le C_{p^n}$ and $V\in RO(H)$.
For
\[
x\in\pi_{q^*V-m\lambda_H}^{\widetilde H}\HZ,
\qquad m\ge0,\quad 2m>\dim V,
\]
write $x/a_{\lambda_H}^m$ for its image
in $\pi_V^H X_{\F_p}$ under
Proposition~\ref{prop:cp-normed-comparison}.
Then
\[
N_H^L\!\left(\frac{x}{a_{\lambda_H}^m}\right)
=
\frac{N_{\widetilde H}^{\widetilde L}(x)}
     {a_{m\Ind_{\widetilde H}^{\widetilde L}\lambda_H}}
\in\pi_{\Ind_H^L V}^L X_{\F_p},
\]
since
\[
\begin{gathered}
N_{\widetilde H}^{\widetilde L}(a_{\lambda_H}^m)
=
a_{m\Ind_{\widetilde H}^{\widetilde L}\lambda_H},
\qquad
\bigl(\Ind_{\widetilde H}^{\widetilde L}\lambda_H\bigr)^{C_p}=0,
\\
\Ind_{\widetilde H}^{\widetilde L}(q^*V)
\cong q^*(\Ind_H^L V).
\end{gathered}
\]

We next determine the norms after extension of coefficients to $A$,
first on pure tensors and then on finite sums.
Under \eqref{cp-flat-extension}, by naturality and \eqref{cp-finite-gauge}, we have
\[
N_H^L(a\otimes x)
=
a^{[L:H]}\otimes N_H^L(x),
\qquad
a\in A,\quad x\in\pi_V^H X_{\F_p}.
\]
For $H=C_{p^h}$ and $L=C_{p^\ell}$, the degree-zero formula
in Witt coordinates is
\[
\begin{aligned}
N_H^L:W_{h+1}(k)&\longrightarrow W_{\ell+1}(k),\\
[\varphi^{-h}(a)]&\longmapsto
[\varphi^{-\ell}(a^{p^{\ell-h}})],
\qquad a\in A,
\end{aligned}
\]
which is the Witt vector norm
\cite[Theorem~4.1]{WittNorm}.

Suppose $[L:H]=p$. For $r\ge2$, translation by $L$ gives
\[
f:L/H\longrightarrow\{1,\ldots,r\},
\qquad |f(L/H)|>1
\quad\Longrightarrow\quad
\operatorname{Stab}_L(f)=H.
\]
Thus \cite[(8.34)]{ZhangHZ} gives $N_H^L$ on finite sums
from products, conjugation and $\tr_H^L$.
For natural $U\in RO(H)$,
\cite[Proposition~6.10]{ZhangHZ},
\eqref{cp-euler-bound} and \eqref{cp-flat-extension} give
\[
\pi_{U-1}^H X_k\cdot\pi_{U-1}^H X_k=0.
\]
Hence
\[
N_H^L(x+y)=N_H^L(x)+N_H^L(y),
\qquad x,y\in\pi_{U-1}^H X_k.
\]

For arbitrary inclusions, use
\[
N_H^L=N_K^L\circ N_H^K
\qquad(H\le K\le L).
\]
\end{proof}

We now compute cyclotomic restriction.
By \cite[Theorems~5.2--5.3]{Ullman}, the identity of $\underline A$ induces an equivalence
\[
d_r:(H\underline A_{C_{p^{r+1}}})^{C_p}
\xrightarrow{\sim}H\underline A_{C_{p^r}},
\qquad r\ge1.
\]
The transition maps of restriction induce
\[
(\mathcal R_{r+1})^{C_p}\xrightarrow{\sim}\mathcal R_r,
\]
which is the identity in Witt coordinates. By
\eqref{cp-witt-gauge}, the map in constant coordinates at
$C_{p^i}$ is
\[
\varphi^{i-1}\circ\operatorname{id}\circ\varphi^{-i}
=\varphi^{-1},\qquad 0\le i\le r.
\]
Naturality of $\omega$ and \cite[Lemma~5.3]{AN} identify
$R_{\mathrm{cyc}}$, under
Proposition~\ref{prop:cp-normed-comparison}, with the map induced
by $H(\varphi^{-1})d_r$.

Let $U=\rho_p(\mathbf a)$, $\mathbf a\in D_n$, and choose an even
$m>0$ with $2m>\dim U+1$. In the complex of
\cite[(5.5), Theorems~5.9, 5.12 and~6.4]{ZhangHZ},
\[
\begin{gathered}
q^*U-m\lambda
=\rho_{C_{p^{n+1}}}(a_0-m,a_0,\ldots,a_n),\\
a_0-m<-\tfrac12,\qquad \ell_0:A_0\xrightarrow{\sim}T_0.
\end{gathered}
\]
Delete $A_0,T_0$ and tensor the remaining complex with $A$ over
$\mathbb Z$. Multiplication by $a_{2\lambda}$ induces
$[q]\mapsto[q]$ on each remaining coordinate
\cite[(A.46)--(A.52)]{ZhangHZ}.
By \eqref{cp-flat-extension} and \eqref{cp-euler-bound},
the kernel and cokernel is $\pi_U^{C_{p^n}}X_k$ and
$\pi_{U-1}^{C_{p^n}}X_k$, respectively.

Assume $n\ge1$. By
\cite[Theorem~4.7 and Corollary~4.35]{ZhangHZ},
the next Euler localization deletes the first remaining
vertex and edge. The target degree is
\[
\varrho_p^*(U^{C_p})
=\rho_{C_{p^{n-1}}}(a_1,\ldots,a_n).
\]
By \cite[(A.46)--(A.52)]{ZhangHZ}, Euler localization induces
$[q]\mapsto[q]$ on each remaining coordinate.
Naturality of the connecting map gives the same projection
on cokernels. For $p=2$,
\[
R_{\mathrm{cyc}}(a_{\sigma_{2^n}})
=\begin{cases}
a_{\sigma_{2^{n-1}}},&n>1,\\
1,&n=1.
\end{cases}
\]
For $p=2$, $n=1$ and $a_1\ge0$, the terminal generator has
fractional representative $2^{a_1}$ in both complexes
\cite[(5.5), (A.52)]{ZhangHZ}.

For $V=V_0+D$, use the sphere unit $\widetilde v_D$ of
Definition~\ref{def:character-translation}.
The unit has nonequivariant degree $1$, so its image
in $H\underline{\mathbb Z}_{(p)}$ is the class of
\cite[Lemma~2.6]{ZhangHZ}.
Inflation commutes with induction along preimage subgroups,
so \cite[(2.7)]{ZhangHZ} applies.
For the maps $f_{r,j}$,
\[
\varrho_p^*(f_{r,j}^{C_p})
=\begin{cases}
\operatorname{id}_{S^0},&j=0,\\
f_{r,j-1},&j>0,
\end{cases}
\]
where $f_{r,j-1}$ on the right is formed over $C_{p^{n-1}}$.
Since $\widetilde v_{r,j}=r f_{r,j}^{-1}$,
\[
\varrho_p^*\Phi^{C_p}(\widetilde v_D)
=\left(\prod_r r^{m_{r,0}}\right)
\widetilde v_{\varrho_p^*(D^{C_p})}.
\]

We now put these results together.

\begin{theorem}\label{thm:cp-result}
For every perfect $\F_p$-algebra $k$ and $n\ge0$, the calculation for
$\HZ$ gives the groups $\pi_V^H X_k$ for all $H\le C_{p^n}$
and $V\in RO(H)$, Tambara functor, and genuine
cyclotomic restriction for $n\ge1$.
\end{theorem}
\begin{proof}
Apply Proposition~\ref{prop:cp-induced} to \cite[Theorems~1.1
and~1.4; (2.7), (5.5), (8.12)--(8.14), (8.28) and (8.31)--(8.34)]{ZhangHZ} gives the tambara functor.
We use the grading convention preceding
Proposition~\ref{prop:vertex-norms}.
For $n\ge1$ and $U=\rho_p(\mathbf a)$, $\mathbf a\in D_n$, put
$W=\varrho_p^*(U^{C_p})$. In degrees
$U\to W$ and $U-1\to W-1$, respectively, $H(\varphi^{-1})d_n$ induces
\[
 R_{\mathrm{cyc}}(x_0,\ldots,x_n)=(\varphi^{-1}x_1,\ldots,\varphi^{-1}x_n),
 \qquad R_{\mathrm{cyc}}[z_0,\ldots,z_{n-1}]=
 [\varphi^{-1}z_1,\ldots,\varphi^{-1}z_{n-1}].
\]
The second target is zero for $n=1$. Multiplicativity and the
fixed-point degrees of $\widetilde v_D$ above extends the map for degree $V_0+D$.
\end{proof}

\enlargethispage{2\baselineskip}
\bibliographystyle{alpha}
\bibliography{references}
\end{document}